\documentclass[11pt,letterpaper]{amsart}

\usepackage{amsmath, amsthm, amssymb, amscd, amsxtra,graphicx}
\usepackage{latexsym, amsfonts}
\usepackage{easybmat}
\usepackage{etex}
\usepackage{url}
\usepackage{texdraw}
\usepackage{epsfig}
\usepackage{tikz}
\allowdisplaybreaks[4]
\usepackage{xcolor, soul}

\usepackage{arydshln}

\usepackage{geometry}
\usepackage{geometry}
\makeatletter \@addtoreset{equation}{section}

\makeatletter \renewcommand{\@biblabel}[1]{#1.}
\theoremstyle{remark}

\begin{document}
\setcounter{page}{1}
\title[One can know the boundary area and curvatures by the Steklov frequencies]{One can know boundary area and curvatures by hearing the Steklov frequencies of a Stokes flow}

\author{Genqian Liu}
\address{School of Mathematics and Statistics, Beijing Institute of Technology, Beijing 100081, China}
\email{liugqz@bit.edu.cn}
\subjclass[2010]{}

\maketitle

\date{}
\protect
\maketitle ~~~\\[-15mm]

\begin{center}
{\footnotesize   School of Mathematics and Statistics, Beijing Institute of Technology, Beijing 100081, China\\
 Emails:  liugqz@bit.edu.cn \\
 }
\end{center}


\vskip 0.48 true cm
\begin{center}
Dedicated to the memory of Professor Louis Nirenberg
\end{center}

\vskip 16.98 true cm

\begin{abstract}
$\,$By calculating full symbol for
   the Dirichlet-to-Neumann map $\Lambda$ of a Stokes flow,
   we establish the asymptotic expansion of the trace of
the heat kernel for $\Lambda$. We also give a useful procedure, by which all coefficients of the asymptotic expansion can be explicitly calculated.
   These coefficients are the Steklov spectral invariants of $\Lambda$, which provide precise geometric information of  the boundary for the Stokes flow. In particular, the first two coefficients show that the area and (total) mean curvature of the the boundary can be known by the Steklov eigenvalues of $\Lambda$.   \end{abstract}

\vskip 1.39 true cm

\section{ Introduction}

\vskip 0.45 true cm

Spectral asymptotics for partial differential operators have been the subject of extensive research for over a century. $\!$It has attracted the attention of many  mathematicians and physicists. Beyond the beautiful asymptotic formulas that are intimately related to the geometric properties of the domain or its boundary, a sustaining force has been its important role in
mathematics, mechanics and theoretical physics (see, for example, \cite{AGMT}, \cite{Cha}, \cite{CLN}, \cite{CH}, \cite{ES}, \cite{Ho4},  \cite{Iv}, \cite{LaST}, \cite{Lo}, \cite{MS}, \cite{Pa1}, \cite{SV}, \cite{Sa}, \cite{Sar}, \cite{SYau}, \cite{Ste}, \cite{We5}).

Let $(\Omega,g)$ be an $n$-dimensional, compact smooth Riemannian manifold with smooth boundary $\partial \Omega$. Assume that $\Omega$ is filled with an incompressible fluid. Let $\mathbf{u} = (u^1,\cdots, u^n)^T$ be the velocity vector field satisfying the stationary Stokes equations
   \begin{eqnarray} \label{2022.3.20-1}  \left\{\!\begin{array}{ll} \mbox{div}\; (\sigma_\mu (\mathbf{u}, p))^\sharp =0 \;\; & \mbox{in} \;\; \Omega,\\
   \mbox{div}\; \mathbf{u}=0 \;\;& \mbox{in}\;\; \Omega,\end{array} \right. \end{eqnarray}
    where $\sigma_{\mu} (\mathbf{u}, p)=2\mu \,(\mbox{Def}\; \mathbf{u})^{\sharp} - p\;\!\mathbf{I}_n$ is the stress tensor,
    $\; \mbox{Def}\; \mathbf{u} = \frac{1}{2}(\nabla \mathbf{u} +\nabla \mathbf{u}^t) $ is the deformation tensor (see Section 2), $\sharp$ is the sharp operator (for a tensor) by raising an index, and $p$ is the pressure. Here $\mu$ is the viscosity function and $\mathbf{I}_n$ denotes the $n\times n$ identity matrix.
Precisely, the first equation of  (\ref{2022.3.20-1}) can be written as (see \cite{Liu5} or Section 2 below)
 \begin{eqnarray} \label{2022.6.24-1} \mu \left( -\nabla^*\nabla u +\mbox{Ric}\, (u) \right) + (\mbox{Def}\; \mathbf{u})^\sharp(\nabla_g \mu)  - \nabla_g p=0 \;\;\, \mbox{in}\;\, \Omega, \end{eqnarray}
   where  $-\nabla^*\nabla \mathbf{u}$  is the Bochner Laplacian of $\mathbf{u}$ (see also \cite{Ta3} or \cite{Liu1}), and $\mbox{Ric}\,(\mathbf{u}) =
   \big( \sum_{l=1}^n R_l^1 u^l$, $\cdots$, $\sum_{l=1}^n R_l^n u^l\big)$ denotes the action of Ricci tensor $R_l^j$ on $\mathbf{u}$.
A fluid flow obeying equations (\ref{2022.3.20-1}) is called the Stokes flow (i.e., creeping flow).
Physically, the viscosity is a function of density, say $\rho$, and temperature.
Here we ignore the effect of temperature. In view of the mass conservation equation, the
incompressibility condition $\mbox{div}\, \mathbf{u}=0$ is equivalent to the fact that the material derivative of
the density function is zero (see, p.$\,$11 of \cite{ChM}), namely,
\begin{eqnarray} \label{2022.6.7-1} \frac{D\rho}{\partial t} := \frac{\partial \rho}{\partial t} +\langle \mathbf{u}, \nabla \rho\rangle =0.\end{eqnarray}
When $\rho$ is a constant, (\ref{2022.6.7-1}) is clearly satisfied. But (\ref{2022.6.7-1}) also holds for nonconstant
density functions. This is indeed the case for spatially stratified fluids, for example, in oceanography. Therefore, for this type of fluid, we expect a spatially varying viscosity $\mu$. We refer the reader to \cite{HLW},  \cite{LUW}, \cite{LiQ}, or  p.$\,$45 of \cite{LL}  for the explanation of viscosity function $\mu$.

Let $\boldsymbol{\phi}\in [H^{\frac{1}{2}} (\partial \Omega)]^n $ satisfy  the standard flux compatibility condition
 \begin{eqnarray} \label{20200502-2} \int_{\partial \Omega} \langle \boldsymbol{\phi},  \boldsymbol{\nu} \rangle \,ds=0 \end{eqnarray}
 where $\boldsymbol{\nu}$ is the unit outer normal field to $\partial \Omega$,  $ds$  denotes the $(n-1)$-dimensional volume element on $\partial \Omega$, and $[H^1(\partial \Omega)]^n= H^1(\partial \Omega) \times \cdots \times H^1(\partial \Omega)$.
 This condition leads to the uniqueness of solution for (\ref{2022.3.20-1}), that is, there exists a unique solution $(\mathbf{u}, p)\in [H^{1} (\Omega)]^n \times H^1(\Omega)$ ($p$ is unique
up to a constant) of (\ref{2022.3.20-1}) and $\mathbf{u}\big|_{\partial \Omega}=\boldsymbol{\phi}$. We will always take $\int_{\Omega} p\, dV=0$ so that $p$ is also unique in the solution of (\ref{2022.3.20-1}) (see \cite{HLW}, \cite{LUW}, \cite{AG} or \cite{Liu5}), where $dV$ denotes the volume element in $\Omega$. It is well-known that the solution of the Stokes equations provides a good approximation to the solution of the Navier-Stokes equations.

$\mathbf{u}\big|_{\partial \Omega}$ is called the Dirichlet boundary condition, and $(\sigma_\mu(\mathbf{u}, p))\boldsymbol{\nu}\big|_{\partial \Omega}: =2
\mu(\mbox{Def}\; \mathbf{u})^\sharp\boldsymbol{\nu} -p\boldsymbol{\nu}$ the Neumann boundary condition. In physics, $\mathbf{u}\big|_{\partial \Omega}$ is the velocity of the Stokes flow on the boundary, and $(\sigma_\mu(\mathbf{u}, p))\boldsymbol{\nu}\big|_{\partial \Omega}$ is the stress acting on $\partial \Omega$ (also called the Cauchy force).
Associated with Stokes flow, the Dirichlet-to-Neumann map $\boldsymbol{\Lambda}$ which
maps $[H^{\frac{1}{2}} (\partial \Omega)]^n $ into $[H^{-\frac{1}{2}} (\partial \Omega)]^n$, is defined by
\begin{eqnarray} \label{2022.3.20-2} \boldsymbol{\Lambda} : \boldsymbol{\phi}\to (\sigma_\mu(\mathbf{u}, p))\boldsymbol{\nu}\big|_{\partial \Omega}, \;\; \, \forall \boldsymbol{\phi} \in [H^{\frac{1}{2}} (\partial \Omega)]^n,\end{eqnarray}
where $(\mathbf{u}, p)$ is the unique solution to (\ref{2022.3.20-1}) satisfying $\mathbf{u}\big|_{\partial \Omega}=\boldsymbol{\phi}$ and $\int_{\Omega} p\, dV=0$.

The Dirichlet-to-Neumann map plays an important role in various inverse boundary value problems. The famous Calder\'{o}n problem associated with Stokes flow asks whether the Dirichlet-to-Neumann map $\Lambda$ uniquely determines the viscosity function $\mu$ in $\Omega$ (see \cite{HLW}, \cite{Liu5}, \cite{IY} or \cite{LUW}).
When $\Omega$ is a real-analytic connected Riemannian manifold with real-analytic boundary, by the result in \cite{Liu5}, the viscosity function $\mu$ is uniquely determined by the Dirichlet-to-Neumann map $\boldsymbol{\Lambda}$.
When $\Omega$ is a bounded domain in flat Euclidean space $\mathbb{R}^n$ with smooth boundary, under an additional assumption, this challenging open problem had been partly answered by R. Lai, G. Uhlmann and J. Wang \cite{LUW} in two dimensional case (see also Imanuvilov and Yamamoto \cite{IY}), and by H. Heck, X. Li, J. Wang \cite{HLW} in three dimensional case. This assumption states that $\mu$ and its normal derivatives on $\partial \Omega$ are uniquely determined by the Ditichlet-to-Neumann map $\boldsymbol{\Lambda}$. The author of present paper in \cite{Liu5} had removed the above additional assumption by discussing the full symbol of another equivalent Dirichlet-to-Neumann map; therefore, the Calder\'{o}n  problem for the Stokes flow has been completely solved in the case of flat Euclidean space \cite{Liu5}. The same result still holds for the Dirichlet-to-Neumann map associated with the Navier-Stokes flow (see \cite{Liu5}).

By Green's formula, it is easy to verify that $\Lambda$ is a self-adjoint, nonnegative, pseudodifferential operator on $[H^{\frac{1}{2}}(\partial \Omega)]^n$. Thus $\Lambda$ has discrete eigenvalues: $0\le \lambda_1 \le \lambda_2 \le \cdots \le \lambda_k \cdots \to +\infty$ with each eigenvalue repeated according to its multiplicity. The corresponding eigenvectors $\{ \boldsymbol{\phi}_j\}$ satisfying $\Lambda \boldsymbol{\phi}_j=\lambda_j \boldsymbol{\phi}_j$, $\; \boldsymbol{\phi}_j\in [C^\infty (\partial \Omega)]^n \cap [H^{\frac{1}{2}} (\partial \Omega)]^n$ can be chosen so that $\{ \mathbf{\phi}_j\}$ forms an orthonormal basis of $[L^2 (\partial \Omega)]^n$.
It is clear that $\lambda_k$  can be characterized variationally as
\begin{align}\label{2022.7.1-2}\\
\!\!\!\!\!\!& \!\!\!\!\!\!\!\lambda_1=\frac{\int_\Omega \left[2\mu \langle \mbox{Def}\, \mathbf{v}_1, \mbox{Def}\, \mathbf{v}_1 \rangle +\mu\, \mbox{Ric}(\mathbf{v}_1, \mathbf{v}_1)+(\mbox{Def}\, \mathbf{v}_1)(\mathbf{v}_1, \nabla_g \mu) \right] dV} {
\int_{\partial \Omega} | \mathbf{v}_1|^2 ds}\nonumber\\
  \!\!\!&\!\!\!\!\!\!\!\!=
  \inf_{\substack{\mathbf{v}\in  [H^1 (\Omega)]^n,\;div\, \mathbf{v}=0\\
  0\ne  \mathbf{v}\in [L^2(\partial \Omega)]^n}} \frac{\int_\Omega \left[2 \mu \langle \mbox{Def}\, \mathbf{v}, \mbox{Def}\, \mathbf{v} \rangle +\mu\, \mbox{Ric}(\mathbf{v}, \mathbf{v})+(\mbox{Def}\, \mathbf{v})(\mathbf{v}, \nabla_g \mu) \right] dV} {
\int_{\partial \Omega} | \mathbf{v}|^2 ds},\;\;\qquad \qquad \nonumber\\
\lambda_k&=\frac{\int_\Omega \left[2\mu \langle \mbox{Def}\, \mathbf{v}_k, \mbox{Def}\, \mathbf{v}_k \rangle +\mu\, \mbox{Ric}(\mathbf{v}_k, \mathbf{v}_k)+(\mbox{Def}\, \mathbf{v_k})(\mathbf{v}_k, \nabla_g \mu) \right] dV} {
\int_{\partial \Omega} | \mathbf{v}_k|^2 ds}\nonumber \\
\!\!\!&\!\!\!\!\!\! \!\!=\!
  \max_{\substack{\mathcal{F}=\{\mathbf{v}\in [H^{1}\! (\Omega)]^n|div\, \mathbf{v}=0\} \\codim(\mathcal {F})=k-1}}
\inf_{\substack{\mathbf{v}\in \mathcal{F}\\ 0\ne  \mathbf{v}\in
[L^2(\partial \Omega)]^n}}\!\!\!\! \!\!\!\!\!\!\!\!\frac{\int_\Omega
\left[ 2\mu \langle \mbox{Def}\, \mathbf{v}, \mbox{Def}\, \mathbf{v} \rangle \!+\!\mu\, \mbox{Ric}(\mathbf{v}, \mathbf{v})\!+\!(\mbox{Def}\, \mathbf{v})(\mathbf{v}, \nabla_{\!g} \mu) \right] dV} { \int_{\partial \Omega}
| \mathbf{v}|^2 ds},\,\, k=2,3, 4,
\cdots \nonumber\end{align}
where $H^m(\Omega)$ is the Sobolev space.

Clearly, the knowledge of the Dirichlet-to-Neumann map associated with the Stokes flow uniquely determines all the eigenvalues $\{\lambda_k\}$  and  the corresponding eigenvectors $\{\boldsymbol{\phi}_k\}$, and vice verse. It is a most fascinating phenomena in nature (a very important physics law) that the ratio of the normal stress acting on the boundary and the corresponding velocity of the fluid on the boundary is a constant if and only if this constant is a Steklov eigenvalue (i.e., a Steklov frequency) for a stationary Stokes flow, and all such ratios (i.e., Steklov frequencies) are discrete. It is not an easy task to simultaneously obtain (or measure) all the eigenvalues and eigenfunctions from the Dirichlet-to-Neumann map associated with Stokes flow. In physics, one can only detect (or measure) all the Steklov frequencies of the Dirichlet-to-Neumann map for a stationary Stokes flow. Thus, the following Kac-type problem for Stokes flow is a quite interesting topic: which geometric quantities can be gotten only by knowing all eigenvalues of the Dirichlet-to-Neumann map of a stationary Stokes flow?

\vskip 0.20 true cm

The main result of this paper is the following:

\vskip 0.19 true cm

\noindent{\bf Theorem 1.1.} \ {\it  Let $\Omega$ be an $n$-dimensional, compact Riemannian manifold with smooth boundary $\partial \Omega$, $n=2,3$.  Assume that $\mu(x)>0$ for all $x\in \bar \Omega$. Let $\{\lambda_k\}$ be the all eigenvalues of the Dirichlet-to-Neumann map $\Lambda$ associated with Stokes flow.
Then \begin{eqnarray} \label{19.6.15-1} \sum_{k=1}^\infty e^{-t\lambda_k} \sim  \sum_{m=0}^{n-1} a_m t^{-n+m+1} + o(1)\, \quad \mbox{as}\;\; t\to 0^+,\end{eqnarray}
where $a_m$ are constants, which can be explicitly calculated by the procedure given in section 5 for $m<n$.
In particular, if $n\ge 2$, then
   \begin{eqnarray} \label{6.0.1}    && \sum_{k=1}^\infty e^{-t\lambda_k} =t^{1-n}
\int_{\partial \Omega} a_0(x') \,ds(x') + t^{2-n}  \int_{\partial \Omega} a_1(x')\, ds(x') \\  &&\quad \quad \quad\quad  \quad + \left\{\!\begin{array}{ll} O(t^{3-n}) \quad \;\, \mbox{when}\;\; n>2,\\ O(t\log t) \quad \, \mbox{when} \;\; n=2,\end{array}\right. \quad \;\mbox{as}\;\; t\to 0^+.\nonumber\end{eqnarray}
Here
   \begin{align}\label{109.6.14-1} &a_0(x')=  \frac{n\, \Gamma(n-1) \,\text{vol}\,( \mathbb{S}^{n-2} ) }{(4\pi )^{n-1} }  \, \mbox{vol}\, (\partial \Omega),\\
\label{06-060} &a_1(x')= \, \frac{(2n+1)\,\Gamma (n-1) \,\text{vol}\,( \mathbb{S}^{n-2} )}{2(n-1)(4\pi )^{n-1}}
\int_{\partial \Omega} \mu(x') \Big(\sum_{\alpha=1}^{n-1}\!{\kappa_{\alpha}(x')}\Big)\, ds(x'),   \end{align}
  where $\kappa_1(x'),\cdots, \kappa_{n-1}(x')$ are the principal curvatures of the boundary $\partial \Omega$ at $x'\in \partial \Omega$, $\text{vol}\left( \mathbb{S}^{n-2} \right) = \frac{2\pi^{(n-1)/2}}{\Gamma(\frac{n-1}{2})} $ is the volume of $(n-2)$-dimensional unit sphere $\mathbb{S}^{n-2}$ in $\mathbb{R}^{n-1}$.}

\vskip 0.35 true cm

Note that the eigenvalues of the  Dirichlet-to-Neumann map $\Lambda$ are physics quantities which can be
measured experimentally. Therefore, Theorem 1.3 shows that the boundary volume ${\mbox{vol}}(\partial \Omega)$, the total mean curvature of the boundary surface $\partial \Omega$ are all spectral invariants and can also be obtained by the Steklov eigenvalues of the fluid.

\vskip 0.25 true cm

The main ideas of this paper are as follows. We first derive a standard expression of the Stokes equations  on a Riemannian manifold in term of a vector field $(\Omega, g)$ (see \cite{Liu5}):
 \begin{eqnarray} \label{2022.6.21-1}\left\{ \!\! \begin{array}{ll} \mu \left( -\nabla^*\nabla \mathbf{u} +\mbox{Ric}\, (\mathbf{u}) \right) +  S^{jk}\,\frac{\partial \mu}{\partial x_k} \frac{\partial }{\partial x_j} - \nabla_g p=0, \;\; &\mbox{in}\;\, \Omega, \\
    \mbox{div}\; \mathbf{u}=0\;\; &\mbox{in} \;\; \Omega,
  \end{array} \right.\end{eqnarray}
  where  $-\nabla^*\nabla \mathbf{u}$ is the Bochner Laplacian of $\mathbf{u}$ (see also \cite{Ta3} or \cite{Liu1}).  Secondly, in boundary normal coordinates, we obtain a  local representation for the Dirichlet-to-Neumann map of the Stokes flow:
  \begin{eqnarray*} \label{2022.6.21-9} \sigma (\mathbf{u}, p) (-\boldsymbol{\nu}) \!\!\! \!\!\!\!\!\! && = 2\mu\, (\mbox{Def}\, \mathbf{u})^\# (-\boldsymbol{\nu}) -p (-\boldsymbol{\nu})\\ [2.5mm]
&&  =\begin{bmatrix}
\begin{BMAT}(@, 15pt, 15pt){c.c}{c.c}
  \! \left[\mu \delta_{jk} \frac{\partial}{\partial x_n} \right]_{(n-1)\times (n-1)}   &  \left[\!\mu g^{j\alpha} \frac{\partial }{\partial x_\alpha}  \!\right]_{(n\!-1)\times \!1} \\
          \left[0 \right]_{1\!\times (n-1)}  & 2\mu \frac{\partial}{\partial x_n}  \end{BMAT}
\end{bmatrix}\!\begin{bmatrix} u^1\\
\cdots \\
u^{n-1} \\
u^n\end{bmatrix} -\begin{bmatrix} 0\\
\cdots \\
0 \\
p\end{bmatrix},\nonumber\end{eqnarray*}
 where $-\boldsymbol{\nu}=(0,\cdots, 0,1)$ is the unit interior normal on $\partial \Omega$.
 The third ingredient is the following transform of functions:
\begin{eqnarray} \label{2022.6.25-7}  \begin{bmatrix} \mathbf{u}\\ p \end{bmatrix} =\begin{bmatrix} (\mu+\rho)^{-\frac{1}{2}} \mathbf{w} +\mu^{-1} \nabla_g f -
 f\nabla_g \mu^{-1} \\ \Delta_g f +\mu ( \Delta_g \mu^{-1}) f + (\mu+\rho)^{-\frac{1}{2}} \frac{\partial \mu}{\partial x_k} w^k  \end{bmatrix},\end{eqnarray}
which transforms (\ref{2022.6.21-1}) into a system of  elliptic equations
    \begin{eqnarray} \label{2022.6.21-10}\;\;\;\;\;\;\;\quad  \left\{ \!\!\begin{array}{ll} L_j(\mathbf{w}, f)=0, \;\; j=1,\cdots,n, \;\; \mbox{in}\;\, \Omega, \\
   \Delta_g f \! -\!\mu (\Delta_g\mu^{-1}) f \!+\! \mu (\mu+\rho)^{-\frac{1}{2}} \frac{\partial w^k }{\partial x_k} \! +
\! \Big( \mu (\mu\!+\!\rho)^{-\frac{1}{2}} \Gamma_{kl}^l\! +\!\mu \frac{\partial ((\mu\!+\!\rho)^{-\frac{1}{2}})}{\partial x_k} \Big) w^k =0\;\,\mbox{in}\;\, \Omega, \end{array} \right. \end{eqnarray}
where $\mathbf{w}=(w^1, \cdots, w^{n})^T$ and $L_j(\mathbf{w},f)$ are given in \cite{Liu5} (see also Section 2).
 By (\ref{2022.6.25-7}), we  can  get the explicit expressions of $\mathbf{u}$ and $(\sigma_\mu (\mathbf{u}, p))(-\boldsymbol{\nu})$ by $(\mathbf{w},f)$:
  \begin{eqnarray*}  \mathbf{u} =\mathbf{L} \begin{bmatrix} \mathbf{w}\\ f \end{bmatrix}, \;\;\,\;\,
    (\sigma_\mu (\mathbf{u}, p))(-\boldsymbol{\nu})= \mathbf{M} \begin{bmatrix} \mathbf{w}\\ f \end{bmatrix}, \end{eqnarray*}
  where $\mathbf{L}$ and $\mathbf{M}$ are $n\times (n+1)$ matrix-valued differential operators (see Section 3).
But, in the matrix-valued differential operators $\mathbf{L}$ and $\mathbf{M}$, there are some terms which contain the forms $\frac{\partial w^1}{\partial x_n}, \cdots, \frac{\partial w^n}{\partial x_n}, \frac{\partial f}{\partial x_n}$. Thus, we need to calculate  the expression of $\frac{\partial (\mathbf{w},f)}{\partial x_n}$. Note that the above elliptic system (\ref{2022.6.21-10}) can define a new Dirichlet-to-Neumann map $\boldsymbol{\Xi}: [H^{\frac{1}{2}}(\partial \Omega)]^{n+1} \to [H^{-\frac{1}{2}}(\partial \Omega)]^{n+1}$ given by
  \begin{eqnarray*} \label{2022.6.21-11} \boldsymbol{\Xi} (\mathbf{w}|_{\partial \Omega}, f|_{\partial \Omega}) = (\frac{\partial \mathbf{w}}{\partial \nu}, \frac{\partial f}{\partial \nu})\big|_{\partial \Omega}, \;\; \; \mbox{for all}\;\, (\mathbf{w}|_{\partial \Omega},f|_{\partial \Omega})\in [H^{\frac{1}{2}}(\partial \Omega)]^{n+1},\end{eqnarray*}
  where $(\mathbf{w},f)$ satisfying (\ref{2022.6.21-10}).
In \cite{Liu5}, by applying the method of factorization (i.e., $ \frac{\partial^2}{\partial x_n^2} \mathbf{I}_{n+1} +\mathbf{B}\frac{\partial }{\partial x_n} +\mathbf{C} = ( \frac{\partial}{\partial x_n} \,\mathbf{I}_{n+1} + \mathbf{B} -\mathbf{Q} ) (\frac{\partial }{\partial x_n}\,\mathbf{I}_{n+1}  +\mathbf{Q})$), the author of the present paper obtained the symbol representation of the pseuddifferential operator $\mathbf{Q}$, that is,
  $\iota(\mathbf{Q})\sim  \sum_{l=0}^\infty \mathbf{q}_{1-l}$, and $\mathbf{q}_{1-l}$, ($l=0,1,2,\cdots$), have been explicitly obtained,  where $\iota(\mathbf{Q})$ denotes the full symbol of $\mathbf{Q}$. It is easy to verify that the pseudodifferential operator $\frac{\boldsymbol{\partial}}{\boldsymbol{\partial} \mathbf{x}_n}:=\mathbf{Q}$ is just the Dirichlet-to-Neumann map $\boldsymbol{\Xi}$.
    Let us point out that unlike classical Dirichlet-to-Neumann map associated with Laplacian, two major difficulties had been encountered for calculating the full symbol $\boldsymbol{\Xi}$. In \cite{Liu5}, these challenging problems had been solved by applying Galois group theory and operator algebra technique, and then by solving several difficult matrix equations (including the famous Sylvester equations).
      Therefore, the key idea in this paper is to seek an $n\times n$ matrix-valued pseudodifferential operator $\boldsymbol{\Psi}$ such that
  \begin{eqnarray*} \boldsymbol{\Psi} \mathbf{L}=\mathbf{ M}. \end{eqnarray*}
 This is equivalent to seek the full symbol  $\iota(\boldsymbol{\Psi}):=\boldsymbol{\psi}$ of a classical pseudodifferential operator $\boldsymbol{\Psi}$ such that $\iota(\boldsymbol{\Psi}\mathbf{L})=\iota(\mathbf{M})$.  Note that $\boldsymbol{\psi}$  has the expansion $\boldsymbol{\psi}\sim \sum_{l=0}^\infty \boldsymbol{\psi}_{1-l}$.
Replacing $\frac{\partial \mathbf{w}}{\partial x_n}$ and $\frac{\partial f}{\partial x_n}$ by the  representation of $-\mathbf{Q}$, and then by solving corresponding matrix equations, we can obtain all $\boldsymbol{\psi}_{1-l}$, $\,l=0,1,2, \cdots$ (see Section 3).
It is easy to check that $\sum_{l=0}^\infty \boldsymbol{\psi}_{1-l}$ is exactly the full symbol of the Dirichlet-to-Neumann map $\boldsymbol{\boldsymbol{\Lambda}}$ associated with Stokes flow.
          Since (see Section 3) \begin{eqnarray*} \boldsymbol{\psi_1}(x',\xi')=2\sqrt{g^{\alpha\beta}(x')\,\xi_\alpha\xi_\beta}\;\,\mathbf{I}_n \end{eqnarray*}
 is a positive definite matrix for all $x'\in \partial \Omega$ and $\xi'\ne 0$, we see that $\boldsymbol{\Lambda}$ is an matrix-valued elliptic pseudodifferential operator of order $1$. This is a surprising result though the Stokes system  is not a system of elliptic equations!
 Therefore, we can consider the heat semigroup $e^{-t \boldsymbol{\Lambda}}$ (associated with $\boldsymbol{\Lambda}$) and its heat trace asymptotic expansion
 as $t\to 0^+$. The coefficients $a_m$ of the asymptotic expansion of the  heat trace are spectral invariants of the operator $\boldsymbol{\Lambda}$ that encode the information about the asymptotic properties of the spectrum. They are of great importance in spectral geometry and have extensive application in physics because they describe true  physical behavior.
The heat invariants have been studied for the Laplacian on Riemannian manifolds with or without boundaries (see \cite{MS}, \cite{Gil} or \cite{Gil2}) as well as for the Dirichlet-to-Neumann map associated with the Laplacian (see \cite{Liu4} and \cite{PoS}) by considering heat
trace asymptotics and applying symbol calculus (see \cite{See2} and \cite{Gr}).
We also refer the reader to \cite{Liu0} and \cite{Liu9} for the asymptotic expansions of the heat traces about the Stokes flow operator and Lame elastic operator in $\Omega$. However, there is no result in the context of the Steklov eigenvalue problem for Stokes flow because of the technical complexity and, most importantly, lack of symbol representation for the ``exotic'' Dirichlet-to-Neumann map $\boldsymbol{\Lambda}$ (In fact, ``most of the studies in spectral geometry and spectral asymptotics are restricted to so-called Laplace type operators'', see p.$\,$120 of \cite{Avr10}). As follows from the results of \cite{Liu4, Liu1, Liu7,GruS}, the trace of the associated heat kernel, $e^{-t \Lambda}$, also admits an asymptotic expansion
\begin{align*} \label{18/12/23-9}  & \sum\limits_{k=0}^\infty e^{-t\tau_k} =\int_{\partial \Omega} \!\bigg\{\frac{1}{(2\pi)^{n-1}} \int_{\mathbb{R}^{n-1}}\!e^{i\langle x'-x',\xi'\rangle}\bigg[\frac{i}{2\pi} \int_{\mathcal{C}} e^{-t \tau}\; \mbox{Tr}\Big( \iota\big((\boldsymbol{\Lambda} -\tau \mathbf{I})^{-1}\big)\Big)d\tau \!\bigg]d\xi'\bigg\} ds(x')\nonumber \qquad \qquad \qquad  \qquad    \\
&\quad\quad\quad  \quad =\int_{\partial \Omega} \!\begin{small}\bigg\{\frac{1}{(2\pi)^{n-1}} \!\int_{\mathbb{R}^{n-1}}\!  \bigg(\!\frac{i}{2\pi}\!\int_{\mathcal{C}} e^{-t \tau} \big( \mbox{Tr}\big({\boldsymbol{\varpi}}_{-1}( x', \xi',\tau)\big)\!+\!\mbox{Tr}
\big( \boldsymbol{\varpi}_{-2} (\tau, x', \xi')\big)+ \cdots \big) d\tau \!\bigg) d\xi'\!\bigg\}\end{small} ds(x')\nonumber\,  \\
&\quad\quad\quad\quad\sim \sum\limits_{m=0}^\infty a_m t^{-n+1+m} +\sum\limits_{l=1}^\infty b_l t^l \log t \quad \mbox{as}\,\, t\to 0^+, \;\quad\nonumber\end{align*}
where $\mathcal{C}$ is a contour around the positive real axis and $\iota((\boldsymbol{\Lambda}-\tau \mathbf{I})^{-1}):=\sum_{m=0}^\infty \boldsymbol{\varpi}_{-1-m}(x',\xi',\tau)$ is the full symbol of pseudodifferential operator $(\boldsymbol{\Lambda} -\tau \mathbf{I})^{-1}$. By using symbol calculus and residue theorem (see Section 5), we can finally get all coefficients $a_m$, $1\le m\le n-1$, for $\boldsymbol{\Lambda}$. These coefficients $a_m$ explicitly give some important geometric information for the Riemannian manifold $\Omega$.

\vskip 0.12 true cm

 The plan of the paper is as follows. In Section
2 we give a standard expression of Stokes equations on a Riemannian manifold and its equivalent system of elliptic equations with $n+1$ unknown functions. In Section 3, by a factorization to this system we obtain the full symbol of a new Dirichlet-to-Neumann map. In Section 4, by solving the corresponding matrix equations, we obtain the full symbol of the Dirichlet-to-Neumann map associated with Stokes flow. Section 5 is devoted to the heat trace asymptotic for the $\boldsymbol{\Lambda}$.

\vskip 1.49 true cm

\section{Stationary Stokes equations on a Riemannian manifold}

\vskip 0.45 true cm

\addvspace{4mm}

We first introduce some concepts of pseudodifferential operators and symbols (see \cite{KN,Ho3,Ho4,Ta2,Gr,Sh}). Assuming $U \subset \mathbb{R}^n$ and $m \in \mathbb{R}$, we define $S^{m}_{1,0} = S^{m}_{1,0}(U,\mathbb{R}^n)$ to consist of $C^{\infty}$-functions $p(x,\xi)$ satisfying for every compact set $W \subset U$,
\begin{align*}
    |D_{x}^{\beta}D_{\xi}^{\alpha}p(x,\xi)|
    \leqslant C_{W,\alpha,\beta} \langle \xi \rangle ^{m-|\alpha|}, \quad x \in W,\ \xi \in \mathbb{R}^n
\end{align*}
for all $\alpha,\beta \in \mathbb{N}^{n}_{+}$, where $D^{\alpha} = D^{\alpha_1} \cdots D^{\alpha_n}$, $D_{j} = \frac{1}{i} \frac{\partial }{\partial x_{j}}$ and $\langle \xi \rangle = (1 + |\xi|^2)^{1/2}$. The elements of $S^{m}_{1,0}$ are called symbols of order $m$. Let $p(x,\xi) \in S^{m}_{1,0}$.  A  pseudodifferential operator in an open set $U$ is essentially defined by a Fourier integral operator
\begin{align*}
    P(x,D)u(x)
    = \frac{1}{(2\pi)^n}\int_{\mathbb{R}^n}p(x,\xi)e^{ix \cdot \xi} \hat{u}(\xi) \,d\xi
\end{align*}
for any $u \in C^{\infty}_{0}(U)$, where $\hat{u}(\xi) = \int_{\mathbb{R}^n}e^{-iy \cdot \xi} u(y) \,dy$ is the Fourier transform of $u$. In such a case we say the associated operator $P(x,D)$ belongs to $OPS^{m}$. If there are smooth $p_{m-j}(x,\xi)$, homogeneous in $\xi$ of degree $m-j$ for $|\xi| \geqslant 1$, that is, $p_{m-j}(x,r\xi) = r^{m-j}p_{m-j}(x,\xi)$, and if
\begin{eqnarray}\label{symbol}
    p(x,\xi) \sim \sum_{j \geqslant 0}p_{m-j}(x,\xi)
\end{eqnarray}
in the sense that
\begin{align*}
    p(x,\xi) - \sum_{j = 0}^{N} p_{m-j}(x,\xi) \in S_{1,0}^{m-N-1}
\end{align*}
for all $N$, then we say $p(x,\xi) \in S_{cl}^{m}$, or just $p(x,\xi) \in S^{m}$. We call $p_{m}(x,\xi)$ the principal symbol of $P(x,D)$. We say $P(x,D) \in OPS^{m}$ is elliptic of order $m$ if on each compact $W \subset U$ there are constants $C_{W}$ and $r < \infty$ such that
\begin{align*}
    |p(x,\xi)^{-1}| \leqslant C_{W} \langle \xi \rangle ^{-m}, \quad |\xi| \geqslant r.
\end{align*}

 An operator $P: C_0^\infty (\Omega) \to C^\infty(\Omega)$ is called a pseudodifferential
operator on a Riemannian manifold $\Omega$ if for any chart diffeomorphism $\varrho: G\to U\subset \mathbb{R}^n$ (where $G\subset \Omega$), the
operator $\tilde P:=  (\varrho^*)^{-1}\circ P \circ\varrho^*$  is a pseudodifferential operator on $U$, where $\varrho^*$ is the pulled back operator of $\varrho$. In other words, it requires $P$ to be locally transformed to pseudodifferential operators on $\mathbb{R}^n$
by some coordinate cover (then P is necessarily so transformed
by every coordinate cover).

  If $A$ and $B$ are two pseudodifferential operators
of order $m$ and $m'$, respectively, then the composition $C=A\circ B$
is a pseudodifferential operator of order $m+m'$
 with the symbol (see, for example,
 (3.17) on p.$\,$13 of \cite{Ta2}, or p.$\,$37 of \cite{Tre})
 \begin{eqnarray} \label{2022.7.1-4}c(x, \xi) \sim \sum_{\vartheta} \frac{i^{|\vartheta|}}{\vartheta!} D^\vartheta_\xi
a(x, \xi)  D^\vartheta_x b(x, \xi),\end{eqnarray} where the sum is taken over all multiindices $\vartheta=(\vartheta_1,\cdots, \vartheta_n)$, $\vartheta!=\vartheta_1!\vartheta_2!\cdots \vartheta_n!$, $|\vartheta|=\vartheta_1+\cdots +\vartheta_n$,  and $a(x, \xi)$ and
$b(x, \xi)$ are the symbols of $A$ and $B$, respectively. In
particular, the principal symbol of $A\circ B$ is
$a_m(x,\xi)b_{m'}(x,\xi)$, where $a_m(x, \xi)$ and $b_{m'}(x, \xi)$  are
the principal symbols of $A$ and $B$, respectively.

\vskip 0.2 true cm

 Throughout this paper, we will use the Einstein summation convention: if the same index name appears exactly twice in any monomial term, once as an upper index and once as a lower index, that term is understood to be summed over all possible values of that index.  We will let Greek indices run from 1 to $n-1$, whereas Roman indices from 1 to $n$, unless otherwise indicated.

Let $\Omega$ be a smooth compact Riemannian manifold of dimension $n$ with smooth boundary $\partial \Omega$.
In the local coordinates $\{x_j\}_{j=1}^n$, we denote by $\bigl\{\frac{\partial}{\partial x_j}\bigr\}_{j=1}^n$ a natural basis for tangent space $T_x \Omega$ at the point $x \in \Omega$.
  Then the Riemannian metric $g$ is given by $g =  g_{jk} \,dx_j \, dx_k$.
  Denote by $[g^{jk}]_{n\times n}$ the inverse of the matrix $[g_{jk}]_{n\times n}$ and set $|g|:= \mbox{det}\, [g_{jk}]_{n\times n}$. In
particular, $d\mbox{V}$, the volume element of $\Omega$ is locally given by $d\mbox{V} = \sqrt{|g|}\, dx_1\cdots dx_n$. By $T\Omega$ and $T^*\Omega$ we denote, respectively, the tangent and cotangent bundle on $\Omega$.  A vector field $\mathbf{X}$ in $T\Omega$  will be denoted as $\mathbf{X}=X^j \frac{\partial }{\partial x_j}$, where $X^j$ is called the $j$-th component of $X$ in given coordinates.
  For smooth vector fields $\mathbf{X}= X^j \frac{\partial}{\partial x_j}\in T\Omega$, $\mathbf{Y}=Y^k
      \frac{\partial}{\partial x_k}\in T\Omega$, the inner product with respect to the metric $g$ is denoted by
\begin{equation*}
    \langle \mathbf{X},\mathbf{Y} \rangle =  g_{jk} X^j Y^k.
\end{equation*}
The divergence operator, in the local coordinates, is denoted by
\begin{equation*}
    \operatorname{div} \mathbf{X}
    = \frac{1}{\sqrt{|g|}}  \frac{\partial}{\partial x_j} \bigl(\sqrt{|g|}\,X^j\bigr),
\end{equation*}
and the gradient operator is denoted by
\begin{equation*}
    \operatorname{grad} u
    =  g^{jk} \frac{\partial u}{\partial x_j} \frac{\partial}{\partial x_k} \quad \text{for}\ u \in C^{\infty}(\Omega),
\end{equation*}
where $(g^{jk}) = (g_{jk})^{-1}$. Thus, one can define the Laplace--Beltrami operator as
\begin{equation*}
    \Delta_{g}
    := \operatorname{div} \operatorname{grad}
    = \frac{1}{\sqrt{|g|}} \frac{\partial}{\partial x_j} \Bigl(\sqrt{|g|}\, g^{jk} \frac{\partial}{\partial x_k}\Bigr).
\end{equation*}

Next, let $\nabla$ be the associated Levi-Civita connection. For each $\mathbf{X} \in T\Omega$, $\nabla \mathbf{X}$ is the tensor of type $(0,2)$ defined by
\begin{eqnarray} \label{9-2.5} (\nabla \mathbf{X})(\mathbf{Y},\mathbf{Z}):= \langle \nabla_{\mathbf{Z}} \mathbf{X}, \mathbf{Y}\rangle, \quad \; \forall\, \mathbf{Y}, \mathbf{Z}\in T\Omega. \end{eqnarray}
It is well-known that in a local coordinate system with the naturally associated frame field on the tangent bundle,
\begin{eqnarray*} \nabla_{\frac{\partial}{\partial x_k}} \mathbf{X} =  \big(\frac{\partial X^j}{\partial x_k} + \Gamma_{lk}^j X^l  \big)\frac{\partial }{\partial x_j}\quad \; \mbox{for}\;\; \mathbf{X}= X^j \frac{\partial}{\partial x_j}, \end{eqnarray*}
 where $\Gamma_{lk}^j= \frac{1}{2} g^{jm} \big( \frac{\partial g_{km}}{\partial x_l} +\frac{\partial g_{lm}}{\partial x_k} -\frac{\partial g_{lk}}{\partial x_m}\big)$ are the Christoffel symbols associated with the metric $g$ (see, for example, \cite{Ta2}). If we denote \begin{eqnarray*} {X^j}_{;k}= \frac{\partial X^j}{\partial x_k} + \Gamma_{lk}^j X^l,\end{eqnarray*}
 then     \begin{eqnarray*} \label{18/10/29}  \nabla_{\mathbf{Y}} \mathbf{X} = Y^k {X^j}_{;k} \,\frac{\partial}{\partial x_j} \;\; \mbox{for}\,\;\mathbf{X}=X^j \frac{\partial }{\partial x_j}, \;\, \mathbf{Y}= Y^{k} \frac{\partial}{\partial x_k}.\end{eqnarray*}
   The symmetric part of $\nabla \mathbf{X}$ is $\mbox{Def}\, \mathbf{X}$, the deformation of $\mathbf{X}$, i.e.,
\begin{eqnarray} (\mbox{Def}\; \mathbf{X})(\mathbf{Y},\mathbf{Z}) =\frac{1}{2} \{ \langle \nabla_{\mathbf{Y}} \mathbf{X}, \mathbf{Z}\rangle +\langle \nabla_{\mathbf{Z}} \mathbf{X}, \mathbf{Y}\rangle \}, \quad \, \forall\, \mathbf{Y}, \mathbf{Z}\in T\Omega\end{eqnarray}
(whereas the antisymmetric part of $\nabla\, \mathbf{X}$ is simply $d\mathbf{X}$, i.e.,
\begin{eqnarray*} d\mathbf{X}(\mathbf{Y}, \mathbf{Z}) = \frac{1}{2} \{ \langle \nabla_{\mathbf{Y}} \mathbf{X}, \mathbf{Z}\rangle -\langle \nabla_{\mathbf{Z}} \mathbf{X}, \mathbf{Y}\rangle \}, \quad \, \forall \,\mathbf{Y}, \mathbf{Z}\in T\Omega.)\end{eqnarray*}
Except for the divergence of a vector field in terms of the
covariant derivative,  one can define a general
notion of divergence of a tensor field (see, p.$\,$148-149 of \cite{Ta1}). If $\Phi$ is a tensor field of type $(k,j)$, with
components $\Phi_{\alpha}^{\beta}=\Phi_{\alpha_1\cdots \alpha_j}\!{}^{\beta_1\cdots \beta_k}$
in a given coordinate system, then $\mbox{div}\, \Phi$ is a tensor field of type $(k-1, j)$, with
components
\begin{eqnarray} \label{2020.8.29-1} \Phi_{\alpha_1\cdots \alpha_j}\!{}^{\beta_1\cdots \beta_{k-1}l}_{\,\,\;\;\;\;\;\;\;\;\;\;\;\;\;\;\;;l}.\end{eqnarray}
 For a scalar function $h$, we also denote $g^{ls} \frac{\partial h}{\partial x_s}$ by $h^{;l}$.
   The Riemann curvature tensor $\mathcal{R}$ of $\Omega$ is given by
\begin{eqnarray} \label{9-2.10} \mathcal{R}(\mathbf{X}, \mathbf{Y})\mathbf{Z} = [\nabla_{\mathbf{X}}, \nabla_{\mathbf{Y}}]\mathbf{Z} -\nabla_{[\mathbf{X},\mathbf{Y}]} \mathbf{Z}, \quad \, \forall\, \mathbf{X}, \mathbf{Y}, \mathbf{Z} \in T\Omega, \end{eqnarray}
where $[\mathbf{X}, \mathbf{Y}] := \mathbf{X}\mathbf{Y} - \mathbf{Y} \mathbf{X}$ is the usual commutator bracket. It is convenient to change this
into a $(0, 4)$-tensor by setting
\begin{eqnarray*} \mathcal{R}(\mathbf{X}, \mathbf{Y}, \mathbf{Z}, \mathbf{W}) := \langle \mathcal{R}(\mathbf{X}, \mathbf{Y})\mathbf{Z}, \mathbf{W}\rangle, \quad \; \forall\, \mathbf{X}, \mathbf{Y}, \mathbf{Z}, \mathbf{W} \in T\Omega.
\end{eqnarray*}
In other words, in a local coordinate system such as that discussed above,
 \begin{eqnarray*} &&R_{jklm}= \Big\langle \mathcal{R}\Big(\frac{\partial}{\partial x_l}, \frac{\partial}{\partial x_m}\Big) \frac{\partial}{\partial x_k}, \frac{\partial}{\partial x_j}\Big\rangle.\end{eqnarray*}
The Ricci curvature $\mbox{Ric}$ on $\Omega$ is a $(0, 2)$-tensor defined as a contraction of $\mathcal{R}$:
\begin{eqnarray*} \mbox{Ric} (\mathbf{X}, \mathbf{Y} ):= \sum_{j=1}^n \Big\langle \mathcal{R}\Big(\frac{\partial }{\partial x_j}, \mathbf{Y}\Big) \mathbf{X}, \frac{\partial }{\partial x_j}\Big\rangle = \sum_{j=1}^n\Big\langle \mathcal{R}\Big(\mathbf{Y},\frac{\partial }{\partial x_j}\Big)\frac{\partial }{\partial x_j}, \mathbf{X}\Big\rangle, \quad \forall\, \mathbf{X}, \mathbf{Y} \in T\Omega.\end{eqnarray*}
  That is,  \begin{eqnarray}\label{19.10.3-2}   R_{jk} =  R^l_{jlk}= g^{lm}R_{ljmk}.\end{eqnarray}

\vskip 0.28 true cm

 On a compact Riemannian manifold, the stationary Stokes equations can be expressed by $1$-form just as the most literatures have done. But we would rather write it in the term of a vector field. The advantage of such an expression is that it is a natural generalization of the Euclidean case and it can be easily understood and calculated in a local frame. Now, assume that the Riemannian manifold $\Omega$ is filled with an incompressible fluid. Let $\mathbf{u} =  u^k \frac{\partial }{\partial x_k}\in T\Omega$ be the velocity vector field.
    Recall that the deformation tensor is a symmetric tensor field of type $(0,2)$ defined by
   \begin{eqnarray*} (\mbox{Def}\; \mathbf{u}) (\mathbf{Y},\mathbf{Z}) = \frac{ 1}{2} \left( \langle \nabla_{\mathbf{Y}} \mathbf{u},  \mathbf{Z}\rangle + \langle \nabla_{\mathbf{Z}} \mathbf{u}, \mathbf{Y}\rangle\right), \,\; \;\; \forall\, \mathbf{Y}, \mathbf{Z}\in T\Omega;\end{eqnarray*}
   in coordinate notation, $(\mbox{Def}\; \mathbf{u})_{jk} =\frac{1}{2} (u_{j;k}+u_{k;j})$, where $u_{j;k}= \frac{\partial u_j}{\partial x_k}-
   \Gamma_{jk}^l u_l$. We have $\mbox{Def}:  C^\infty (\bar \Omega, T) \to C^\infty (\bar \Omega, S^2 T^* )$ (see p.$\,$464 of \cite{Ta1}). This tensor was introduced in Chap.$\,$2, $\S 3$, cf (3.35) of \cite{Ta1}.  It follows from p.$\,$153 of \cite{Ta1} (or p.$\,$305 of \cite{Ta2}) that   \begin{eqnarray} \label{20200507-1}  \mbox{Def}^*w=-\mbox{div}\, (w^{\sharp\sharp})
\end{eqnarray}
    for any (0,2) type tensor $w$, where  $\sharp\sharp$ is the twice sharp operator (Note that in local coordinates, $w^{\sharp\sharp} = w^{jk}=g^{jl}g^{km} w_{lm}$ for any tensor $w=w_{jk}$ of type (0,2)). In other words, the adjoint  ${\mbox{Def}}^*$ of $\mbox{Def}$ is defined in local coordinates by \begin{eqnarray} \label{20200507-2} ({\mbox{Def}}^* {w})^j=- {{w}}^{jk}_{\;\;\;\,;k}\end{eqnarray}  for each symmetric tensor field ${w}:=w_{jk}$ of type $(0,2)$.
 In particular, if $\boldsymbol{\nu}\in T\Omega$ is the outward unit normal to $\partial \Omega\hookrightarrow \Omega$, then the integration by parts formula (see formula (2.16) of \cite{DMM}, or formula (12.4) of p. 463 in \cite{Ta1})
\begin{eqnarray} \label{18/113-4} \int_{\Omega} \langle \mbox{Def}\; \mathbf{u}, w\rangle dV= \int_\Omega \langle \mathbf{u}, \mbox{Def}^* w\rangle dV+ \int_{\partial \Omega} w( \boldsymbol{\nu}, \mathbf{u}) \,ds\end{eqnarray}
holds for any $\mathbf{u}\in T\Omega$ and any symmetric tensor field $w$ of type $(0,2)$.
   Thus we have the following:

 \vskip 0.29 true cm

 \noindent{\bf Theorem 2.1} (see \cite{Liu5}).  {\it On a Riemannian manifold $\Omega$, modeling a nonhomogeneous, linear, incompressible fluid, the stationary Stokes equations is given by}
 \begin{eqnarray} \label{2020.4.18-1}\left\{\!\!\begin{array}{ll}  \mbox{div}\;(\sigma_{\mu} (\mathbf{u},p))^\sharp=0 \quad &\mbox{in}\;\; \Omega,\\
  \mbox{div}\; \mathbf{u}=0\quad &\mbox{in} \;\; \Omega,
  \end{array} \right.\end{eqnarray}
   {\it where} $\sigma_{\mu} (\mathbf{u}, p)=2\mu \,(\mbox{Def}\, (\mathbf{u}))^{\sharp} - p\mathbf{\,I}_n$,
    $(\mbox{Def}\;\mathbf{u})^\sharp$ {\it is a tensor of field of type $(1,1)$}. {\it Or equivalently, (\ref{2020.4.18-1}) can be written as}
   \begin{eqnarray} \label{2020.8.31-4}\left\{ \!\! \begin{array}{ll} \mu \left( -\nabla^*\nabla \mathbf{u} +\mbox{Ric}\, (\mathbf{u}) \right) +  S^{jk}\,\frac{\partial \mu}{\partial x_k} \frac{\partial }{\partial x_j} - \nabla_g p=0, \;\; &\mbox{in}\;\, \Omega, \\
    \mbox{div}\; \mathbf{u}=0\;\; &\mbox{in} \;\; \Omega,
  \end{array} \right.\end{eqnarray}
 {\it where}  $-\nabla^*\nabla \mathbf{u}$ {\it is the Bochner Laplacian of $\mathbf{u}$ defined by (see \cite{Liu5}, \cite{Ta3} or \cite{Liu1})}  \begin{eqnarray}\label{2020.4.18-5} && -\nabla^*\nabla \mathbf{u} = \left\{ \Delta_g u^j +2  g^{kl} \Gamma_{sk}^j \frac{\partial u^s}{\partial x_l}
 + \Big( g^{kl} \frac{\partial \Gamma_{sl}^j}{\partial x_k} + g^{kl} \Gamma_{hl}^j \Gamma_{sk}^h -
g^{kl}\Gamma_{sh}^j \Gamma_{kl}^h \Big) u^s\right\} \frac{\partial}{\partial x_j},\end{eqnarray}
    and  \begin{eqnarray} \label{18/11/25} \mbox{Ric} (\mathbf{u})= R_l^{\,j} u^l\,\frac{\partial }{\partial x_j},  \;\;\mbox{and}\;\; S^{jk}=\frac{1}{2}\big(u^{j;k}+u^{k;j}\big). \end{eqnarray}
  {\it In particular, for the velocity vector field  $\mathbf{u}\in T\Omega$, natural boundary conditions for the stationary Stokes equations include prescribing $\mathbf{u}\big|_{\partial \Omega}$, Dirichlet type, and}
  \begin{eqnarray}\label{18/11/6;5}  (\sigma_\mu (\mathbf{u},p))\boldsymbol{\nu} := 2\mu \,(\mbox{Def}\;\mathbf{u})^\sharp \boldsymbol{\nu} - p\boldsymbol{\nu} \quad \mbox{on}\; \, \partial \Omega, \end{eqnarray}
  {\it Neumann type}.

\vskip 0.25 true cm

 For a Riemannian manifold $\Omega$, let $ \boldsymbol\phi\in [ H^{\frac{1}{2}} (\partial \Omega) ]^n$ satisfy $\int_{\partial \Omega} \langle \boldsymbol{\phi}, \boldsymbol{\nu}\rangle ds =0$, then there exists  a unique $(\mathbf{u},p) \in [H^2(\Omega)]^n \times H^1(\Omega)$  solving (\ref{2020.4.18-1}) (see, for example, \cite{Liu5}, \cite{Ter} or A of Chapter 17 in \cite{Ta3}) and
 \begin{eqnarray} \label{2020.6.26-5} \mathbf{u}\big|_{\partial \Omega} =\boldsymbol{\phi}, \,\; \int_\Omega p\, dV=0.\end{eqnarray} So we can naturally define the Dirichlet-to-Neumann map $\Lambda$ of the Stokes flow by
\begin{eqnarray} \label{2022.6.26-2} \boldsymbol{\Lambda} : \boldsymbol{\phi}\to (\sigma_\mu(\mathbf{u}, p))\boldsymbol{\nu}\big|_{\partial \Omega}, \;\; \, \forall \boldsymbol{\phi} \in [H^{\frac{1}{2}} (\partial \Omega)]^n,\end{eqnarray} where $(\mathbf{u},p)$ satisfies (\ref{2020.4.18-1}) and (\ref{2020.6.26-5}).

 For any $\mathbf{v}=(v^1, \cdots, v^n)^T\in [H^1(\Omega)]^n$ satisfying $\mbox{div}\,\mathbf{v}=0$ in $\Omega$,  by taking inner with $\mathbf{v}$ in (\ref{2020.8.31-4}), and then by applying Green's formula we get (see, for example, p.$\,$586 of  \cite{Ta3})
 \begin{eqnarray*} \int_\Omega 2\mu \langle \mbox{Def}\; \mathbf{u}, \mbox{Def}\; \mathbf{v}\rangle dV- \int_{\partial \Omega}\langle \sigma_\mu (\mathbf{u},p), \mathbf{v}\rangle ds + \int_{\Omega} \mu u^j R_{jk} v^k dV + \int_{\Omega}( \nabla_g \mu )^jS_{jk} v^k\,dV=0.\end{eqnarray*}
The above equality can be re-written as
\begin{eqnarray}\label{2022.7.1-1}&& \int_{\partial \Omega}\langle \sigma_\mu (\mathbf{u},p), \mathbf{v}\rangle ds= \int_\Omega 2\mu \langle \mbox{Def}\; \mathbf{u}, \mbox{Def}\; \mathbf{v}\rangle dV + \int_{\Omega} \mbox{Ric}\,(\mathbf{u}, \mathbf{v}) dV \\
&&\;\; \;\;\;+ \int_{\Omega}(\mbox{Def}\; \mathbf{u})( \nabla_g \mu, \mathbf{v} )\,dV \;\; \;\mbox{for any}\,\; \mathbf{v}\in [H^{1}(\Omega)]^n\cap\{\tilde{\mathbf{v}}\in [H^1(\Omega)]^n\big|\mbox{div}\, \tilde{\mathbf{v}}=0\}.\nonumber\end{eqnarray}
From this, we can immediately get the variational expressions (\ref{2022.7.1-2}) for the Steklov eigenvalues $\{\lambda_k\}_{k=1}^\infty$ of  $\boldsymbol{\Lambda}$.

\vskip 0.29 true cm

Next, as be shown in \cite{Liu5}, we use the transform of functions:
      \begin{eqnarray} \label{200420-12}  \begin{pmatrix} \mathbf{u}\\ p \end{pmatrix} =\begin{pmatrix} (\mu+\rho)^{-\frac{1}{2}} \mathbf{w} +\mu^{-1} \nabla_g f -
 f\nabla_g \mu^{-1} \\ \Delta_g f +\mu ( \Delta_g \mu^{-1}) f + (\mu+\rho)^{-\frac{1}{2}} \frac{\partial \mu}{\partial x_k} w^k  \end{pmatrix}.\end{eqnarray}
It follows (see \cite{Liu5}) that $(\mathbf{u},p)$ is a solution of the stationary Stokes equations (\ref{2020.4.18-1}) provided   $(\mathbf{w},f)$
satisfies
    \begin{eqnarray} \label{200421-13}\;\;\;\;\;\;\;\;\;\; \left\{ \!\!\begin{array}{ll} L_j(\mathbf{w}, f)=0, \,\quad j=1,\cdots,n, \;\, \mbox{in}\;\, \Omega,\\
   \Delta_g f \! -\!\mu (\Delta_g\mu^{-1}) f \!+\! \mu (\mu+\rho)^{-\frac{1}{2}} \frac{\partial w^k }{\partial x_k} \! +
\! \Big( \mu (\mu\!+\!\rho)^{-\frac{1}{2}} \Gamma_{kl}^l\! +\!\mu \frac{\partial ((\mu\!+\!\rho)^{-\frac{1}{2}})}{\partial x_k} \Big) w^k =0\;\; \mbox{in}\;\, \Omega, \end{array} \right. \end{eqnarray}
where
   \begin{eqnarray}    &&   \label{200422-3}    L_j(\mathbf{w},f):    = \Delta_g w^j +\bigg[ \mu (\mu^{-1})^{;j} \, \frac{\partial w^k}{\partial x_k}  + \frac{\rho}{\mu(\mu+\rho)} \frac{\partial \mu}{\partial x_m} g^{ml} \frac{\partial w^j}{\partial x_l} +2 g^{ml} \Gamma_{km}^j \frac{\partial w^k}{\partial x_l}\bigg]  \\
   &&   + \bigg[ \Big(  (\mu+\rho)^{\frac{1}{2}} \Delta_g ((\mu+\rho)^{-\frac{1}{2}} ) +\mu^{-1} (\mu+\rho)^{\frac{1}{2}}  \frac{\partial \mu}{\partial x_l} \big((\mu+\rho)^{-\frac{1}{2}}\big)^{;l} \Big) w^j +  \mu (\mu^{-1})^{;j}  \Gamma_{kl}^l w^k \bigg.\nonumber\\
&& \bigg.+ \mu (\mu+\rho)^{\frac{1}{2}} \,  (\mu^{-1})^{;j} \,\frac{\partial ((\mu+\rho)^{-\frac{1}{2}})}{\partial x_k} w^k
+ \frac{\rho}{\mu(\mu+\rho)} g^{ml} \Gamma_{km}^j \frac{\partial \mu}{\partial x_l} w^k\bigg.
 \nonumber
\\
&&\bigg. +   \Big( g^{ml} \frac{\partial \Gamma_{kl}^j}{\partial x_m} + g^{ml} \Gamma_{hl}^j \Gamma_{km}^h -
g^{ml}\Gamma_{kh}^j \Gamma_{ml}^h \Big)  w^k  +R^j_k \, w^k  - \mu^{-1}  \big(\frac{\partial \mu}{\partial x_k}\big)^{;j} w^k \bigg]\nonumber \\
&& +2\mu^{-1} (\mu+\rho)^{\frac{1}{2}}\Gamma_{sr}^j g^{sl}g^{rm}\frac{\partial^2 f}{\partial x_l\partial x_m} + \bigg[ 2\mu^{-1} (\mu+\rho)^{\frac{1}{2}}  R^j_m g^{lm}\frac{\partial f}{\partial x_l}  \bigg. \nonumber\\
    && \bigg. -2(\mu+\rho)^{\frac{1}{2}} g^{ml}  \frac{\partial ((\mu^{-1})^{;j})}{\partial x_m} \, \frac{\partial f}{\partial x_l}  -2(\mu+\rho)^{\frac{1}{2}} g^{ml} \Gamma_{sm}^j (\mu^{-1})^{;s} \frac{\partial f}{\partial x_l}\bigg.
 \nonumber\\
&&\bigg. +   \mu^{-1} (\mu+\rho)^{\frac{1}{2} } \Big( g^{mr} \frac{\partial \Gamma_{sr}^j}{\partial x_m} - g^{mr} \Gamma_{hr}^j \Gamma_{sm}^h  -
g^{mr}\Gamma_{sh}^j \Gamma_{mr}^h \Big) g^{sl} \frac{\partial f}{\partial x_l}\bigg.\nonumber
\\ &&\bigg.
 -2\mu^{-1} (\mu+\rho)^{\frac{1}{2}}
\Gamma_{sh}^j g^{sr}g^{hm}\Gamma_{rm}^l \, \frac{\partial f}{\partial x_l}
 \bigg] +\bigg[ - 2  (\mu+\rho)^{\frac{1}{2}}  (\Delta_g \mu^{-1})^{;j} f \bigg. \nonumber \\
&& \bigg. - 2(\mu+\rho)^{\frac{1}{2}} R^j_l\, (\mu^{-1})^{;l} f
-   2(\mu+\rho)^{\frac{1}{2}} g^{ml}\Gamma_{sm}^j\, \frac{\partial ((\mu^{-1})^{;s}\!)}{\partial x_l}\, f  - (\mu+\rho)^{\frac{1}{2} } \Big( g^{ml} \frac{\partial \Gamma_{sl}^j}{\partial x_m} \bigg.\nonumber\\
&&\bigg. + g^{ml} \Gamma_{hl}^j \Gamma_{sm}^h -
g^{ml}\Gamma_{sh}^j \Gamma_{ml}^h \Big)(\mu^{-1})^{;s} f
  -   2 \mu^{-1} (\mu+\rho)^{\frac{1}{2} }
   \frac{\partial \mu}{\partial x_m} \, (\mu^{-1})^{;m;j} f \bigg].
   \nonumber
       \end{eqnarray}
Clearly, (\ref{200421-13}) is a system of second-order  linear elliptic equations with unknown functions $w^1$, $\cdots$,
$w^n$, $f$ in $\Omega$.
 We further consider the following two elliptic boundary value problems:
  \begin{eqnarray}\label{20200524-3} \;\;\;\;\; \;\;\;\,\left\{\!\!\! \begin{array}{ll} L_j(\mathbf{w}, f)=0, \quad j=1,\cdots,n \quad &\mbox{in}\;\, \Omega, \\
     \Delta_g f \! -\!\mu (\Delta_g\mu^{\!-1}) f \!+\! \mu (\mu\!+\!\rho)^{-\frac{1}{2}} \frac{\partial w^k }{\partial x_k} \! +\!
 \Big(\! \mu (\mu\!+\!\rho)^{\!-\frac{1}{2}} \Gamma_{kl}^l \!+\!\mu \frac{\partial ((\mu\!+\!\rho)^{\!-\frac{1}{2}})}{\partial x_k}\! \Big) w^k\!=\!0 \;\, &\mbox{in}\;\; \Omega,\\
 (\mu+\rho)^{-\frac{1}{2}} w +\mu^{-1} \nabla_g f -
 f\nabla_g \mu^{-1}=u_0 \;\, &\mbox{on}\;\, \partial \Omega,\\
-\mu (\mu\!+\!\rho)^{-\frac{1}{2}} \frac{\partial w^j}{\partial x_j}\! +\!2\mu (\Delta_g \mu^{-1})f \!+\! (\mu \!+\!\rho)^{-\frac{1}{2}}\Big(\!\!-\!\mu \Gamma^j_{kj} \!+\!\frac{3\mu\!+\!2\rho}{2(\mu\!+\!\rho)} \, \frac{\partial \mu}{\partial x_k} \!\Big)w^k\!=\!p_0 \;\; &\mbox{on}\;\, \partial \Omega
  \end{array} \right.
  \end{eqnarray}
 and
  \begin{eqnarray}\label{20200524-4}\; \;\quad\;\; \quad  \left\{\!\!\! \begin{array}{ll} L_j(\mathbf{w}, f)=0, \quad j=1,\cdots,n \quad &\mbox{in}\;\; \Omega, \\
     \Delta_g f \! -\!\mu (\Delta_g\mu^{\!-1}) f \!+\! \mu (\mu\!+\!\rho)^{-\frac{1}{2}} \frac{\partial w^k }{\partial x_k} \! +\!
 \Big(\! \mu (\mu\!+\!\rho)^{\!-\frac{1}{2}} \Gamma_{kl}^l \!+\!\mu \frac{\partial ((\mu\!+\!\rho)^{\!-\frac{1}{2}})}{\partial x_k}\! \Big) w^k\!=\!0 \;\; &\mbox{in}\;\; \Omega,\\
 (\mathbf{w},f)= (\mathbf{w}_0,f_0) \;\; &\mbox{on}\;\, \partial \Omega. \end{array} \right.
  \end{eqnarray}
 Denote by $\mathcal{P}$ the differential operator (defined in $\Omega$) on the left-hand side of  (\ref{20200524-3}) (or (\ref{20200524-4})).
 If we discuss the spectrum of the operator $\mathcal{P}$ in $\Omega$ with vanishing boundary conditions corresponding to the above two systems, we see that the spectrum of $\mathcal{P}$ are all discrete eigenvalues, and any eigenvalue of each kind problem will continuously vary in $\rho$.
   Thus we can choose a suitable constant $\rho\ge 0$ such that $\tau=0$ is not an eigenvalue of the following two eigenvalue problems:
    \begin{eqnarray}\label{20200627-4} \;\; \;\;\;\;\;\quad \left\{\!\!\! \begin{array}{ll} \mathcal{P}(\mathbf{w}, f)= -\tau (\mathbf{w},f) \;\; &\mbox{in}\;\; \Omega,\\
 (\mu+{\rho})^{-\frac{1}{2}} \mathbf{w} +\mu^{-1} \nabla_g f -
 f\nabla_g \mu^{-1}=0 \;\; &\mbox{on}\;\, \partial \Omega,\\
-\mu (\mu\!+\!{\rho})^{-\frac{1}{2}} \frac{\partial w^j}{\partial x_j}\! +\!2\mu (\Delta_g \mu^{-1})f \!+\! (\mu \!+\!\rho)^{-\frac{1}{2}}\Big(\!\!-\!\mu \Gamma^j_{kj}\! +\!\frac{3\mu\!+\!2{\rho}}{2(\mu\!+\!{\rho})} \, \frac{\partial \mu}{\partial x_k} \!\Big)w^k=0 \;\; &\mbox{on}\;\, \partial \Omega
  \end{array} \right.
  \end{eqnarray}
 and
  \begin{eqnarray}\label{20200627-5}  \left\{\!\!\! \begin{array}{ll} \mathcal{P}(\mathbf{w}, f)= -\tau (\mathbf{w},f) \;\; &\mbox{in}\;\; \Omega,\\
 (\mathbf{w},f)= 0 \;\; &\mbox{on}\;\, \partial \Omega. \end{array} \right.
  \end{eqnarray}
  It follows that for any $(\mathbf{w}_0,f_0)\in [H^{\frac{1}{2}} (\partial \Omega)]^{n}\times H^{\frac{1}{2}}(\partial \Omega)$, there is a unique solution $(\mathbf{w},f)\in [H^{1}(\Omega)]^n \times H^1(\Omega)$ of the system (\ref{20200524-4}) satisfying $(\mathbf{w},f)\big|_{\partial \Omega}=(\mathbf{w}_0,f_0)$. Thus, we can define a new Dirichlet-to-Neumann map ${\boldsymbol{\Xi}}: [H^{\frac{1}{2}} (\partial \Omega)]^{n}\times H^{\frac{1}{2}}(\partial \Omega)\to [H^{-\frac{1}{2}} (\partial \Omega)]^{n}\times H^{-\frac{1}{2}}(\partial \Omega)$  associated with
  system (\ref{20200524-4}) by
\begin{eqnarray} {\boldsymbol{\Xi}} (\mathbf{w}_0,f_0) =\frac{\partial (\mathbf{w}, f)}{\partial \nu}\big|_{\partial \Omega}  \quad \; \mbox{for any}\;\; (\mathbf{w}_0,f_0)\in [H^{\frac{1}{2}}(\partial \Omega)]^n\times H^{\frac{1}{2}}(\partial \Omega),\end{eqnarray}
where $(\mathbf{w},f)$ satisfies (\ref{20200524-4}).

\vskip 1.49 true cm

\section{Factorization of new system}

\vskip 0.45 true cm

   From now on, we will denote by
\begin{equation*}
  \begin{bmatrix}
\begin{BMAT}(@, 15pt, 15pt){c.c}{c.c}
   [a_{jk}]_{n\times n}&[b_j]_{n\times 1}\\
   [c_k]_{1\times n}& d
  \end{BMAT}
\end{bmatrix}  \;\;\;\mbox{and}\;\; \left[\!
\begin{array}{c;{3pt/3pt}c} \left[ a_{jk}\right]_{n\times n} & \left[b_j\right]_{n\times 1}\end{array} \!\right]
\end{equation*}
the block matrix
\[ \begin{bmatrix}
\begin{BMAT}(@, 15pt, 15pt){ccc.c}{ccc.c}
    a_{11} & \cdots & a_{1n} & b_{1}\\
 \cdots & \cdots &\cdots & \vdots \\
   a_{n1} & \cdots & a_{nn} & b_{n}\\
    c_{1} & \cdots & c_{n} & d
  \end{BMAT}
\end{bmatrix} \;\;\; \mbox{and}\;\; \left[\!
\begin{array}{ccc;{3pt/3pt}c} a_{11} & \cdots & a_{1n} & b_1 \\
\vdots  & & \vdots & \vdots  \\
 a_{n1} & \cdots & a_{nn} & b_n \end{array}\right] \]
 where $\big[a_{jk}\big]_{n\times n}$, $\big[ b_j\big]_{n\times 1}$ and $\big[c_k\big]_{1\times n}$ are the $n\times n$ matrix \begin{eqnarray*}\begin{bmatrix} a_{11}& a_{12} & \cdots &a_{1n}\\
a_{21}& a_{22} & \cdots &a_{2n}\\
\cdots & \cdots& \cdots & \cdots\\
a_{n1}& a_{n2} & \cdots &a_{nn}\end{bmatrix}, \end{eqnarray*}
the $n\times 1$ matrix
\begin{eqnarray*}\begin{bmatrix} b_1 \\
\vdots \\ b_n\end{bmatrix} \end{eqnarray*}
and the $1\times n$ matrix \begin{eqnarray*}\begin{bmatrix} c_1 & \cdots & c_{n} \end{bmatrix}, \end{eqnarray*} respectively.

In the local coordinates,  we can rewrite (\ref{200421-13}) as
\begin{eqnarray} \label{2022.6.27-10}\end{eqnarray}
\begin{align*}
 & \left\{\begin{bmatrix}
\begin{BMAT}(@, 15pt, 15pt){c.c}{c.c}
   [\delta_{jk}\Delta_g ]_{n\times n}&\left[2\mu^{-1} (\mu+\rho)^{\frac{1}{2}}\Gamma_{sr}^j g^{sl}g^{rm} \frac{\partial^2}{ \partial x_l \partial x_m}\right]_{n\times 1}\\
   [0]_{1\times n}& \Delta_g
  \end{BMAT}
\end{bmatrix}\right.
\\
&   \!\! +\!\begin{small} \begin{bmatrix}\!
\begin{BMAT}(@, 1pt, 5pt){c.c}{c.c}
   \!\left[\! \!\! \begin{array}{ll}  \mu (\mu^{\!-\!1}\!)^{;j}\!\frac{\partial }{\partial x_k}\!+\! \frac{\rho \delta_{jk}}{\mu(\mu\!+\!\rho)} \frac{\partial \mu}{\partial x_m} g^{ml} \frac{\partial}{\partial x_l}  \\
   + 2 g^{ml}\Gamma_{\!km}^j\frac{\partial }{\partial x_l} \end{array} \!\right]_{n\times n}   &   \! \left[\!\frac{(\mu\!+\!\rho)^{\!\frac{1}{2}}}{\mu}\!\! \left\{\!\!\! \begin{array} {ll} \! 2R^j_m g^{lm} \!-\!2\mu g^{ml}\frac{\partial ((\mu^{\!-\!1})^{;j}\!)}{\partial x_m} \!  -\!2 \mu g^{ml} \Gamma_{\!sm}^j (\mu^{\!-\!1})^{;s}\\
        \!+ \!( g^{mr}\frac{\partial \Gamma_{sr}^j }{\partial x_m}  \!-\!  g^{mr} \Gamma_{\!hr}^j \Gamma_{\!sm}^h \!-\! g^{mr} \Gamma_{\!sh}^j \Gamma_{\!mr}^h) g^{sl}  \\
       \! -\!2\Gamma_{sh}^j g^{sr}g^{hm}\Gamma_{rm}^l\end{array}\!\!\!\! \!\right\}\!\!\frac{\partial}{\partial x_l}\!\right]_{\!n\!\times\! 1}\\
   \left[\mu(\mu+\rho)^{\!-\frac{1}{2}} \frac{\partial}{\partial x_k} \!\right]_{1\times n}& 0
  \end{BMAT}\!
\end{bmatrix}\end{small}
\\
&    +\begin{bmatrix}
\begin{BMAT}(@, 5pt, 5pt){c.c}{c.c}
   \!\left[\!\Big((\mu\!+\!\rho)^{\!\frac{1}{2}} \Delta_g ((\mu+\rho)^{\!-\frac{1}{2}}) +\mu^{\!-\!1} (\mu+\rho)^{\frac{1}{2}}   \frac{\partial \mu}{\partial x_l}   ((\mu+\rho)^{\!-\frac{1}{2}})^{;l}\Big) \delta_{jk}  \right]_{n\times n}   &   \! \left[\, 0\,\right]_{n\times 1}\\
   \left[0 \right]_{1\times n}&  -\mu (\Delta_g \mu^{-1})
  \end{BMAT}
\end{bmatrix}
\\
\!\!& \! \begin{small}\left. \! +\!\begin{bmatrix}
\begin{BMAT}(@, 1pt, 1pt){c.c}{c.c}
  \! \!\left[\!\!\!\left. \begin{array}{ll}\! \mu(\mu^{\!-\!1}\!)^{;j} \Gamma_{\!kl}^l \!+\!\mu(\mu\!+\!\rho)^{\!\frac{1}{2}}(\mu^{\!-\!1})^{;j} \frac{\partial (\mu\!+\!\rho)^{\!-\!\frac{1}{2}}}{\partial x_k}\\
 \! \!+\! \big( g^{ml} \frac{\partial \Gamma_{\!kl}^j}{\partial x_m}\! + \!g^{ml}\Gamma_{\!hl}^j \Gamma_{\!km}^h \!\!-\! g^{ml} \Gamma_{\!kh}^j \Gamma_{\!ml}^h \big)\\
   \frac{\rho}{\mu(\mu\!+\!\rho)} g^{ml} \Gamma_{\!km}^j \frac{\partial \mu}{\partial x_l} \!+\!R^j_k \!-\! \mu^{\!
-\!1}\big( \!\frac{\partial \mu}{\partial x_k}\big)^{;j}\end{array} \!\!
\right.\!\! \! \right]_{\!n\!\times \!n}   &   \! \left[\!(\mu\!+\!\rho)^{\!\frac{1}{2}}\!\!\left\{\!\!\! \begin{array}{ll} \!\!-2 ( \Delta_g \mu^{\!-\!1})^{;j} - 2  R_l^j \,(  \mu^{\!-1})^{;l}  \\ \!\!-2  g^{ml} \Gamma_{\!sm}^j \frac{\partial((\mu^{\!-\!1})^{;s})}{\partial x_l}
  \\ -2 \mu^{\!-\!1}\frac{\partial \mu}{\partial x_m} (\mu^{\!-\!1}\!)^{;m;j}
    \!- \!\big(\! g^{ml} \frac{\partial \Gamma^j_{\!sl}}{\partial x_m}\\
     + \!g^{ml} \Gamma^j_{\!hl}\Gamma^h_{\!sm} \!\!-\! g^{ml} \Gamma_{\!sh}^j\Gamma^h_{\!ml}\! \big) \!(u^{\!-\!1}\!)^{;s}\end{array}\!\!\!\!\!\right\}\!\!\right]_{\!n\!\times \!1}\\
   \left[\mu(\mu+\rho)^{-\frac{1}{2}} \Gamma_{kl}^l + \mu \frac{\partial ((\mu+\rho)^{-\frac{1}{2}})}{\partial x_k} \! \right]_{\!1\times n}&  0
  \end{BMAT}
\!\end{bmatrix}\!\right\}\end{small}\!\!\begin{bmatrix} w^1\\
\vdots\\
w^n\\
f\end{bmatrix}\!\!=\!0,
\end{align*}
 where  \[ \delta_{jk}= \left\{\begin{array}{ll} 1 \;\; &\mbox{for}\;\, j=k\\
 0 \;\; &\mbox{for} \,\; j\ne k \end{array} \right.\]
 is the standard Kronecker symbol.

In order to describe the Dirichlet-to-Neumann map associated with  the equivalent new system, we first recall the construction of usual geodesic coordinates with respect to the boundary (see p.$\,$1101 of \cite{LU}, \cite{Liu3} or \cite{Liu4}). For each $x'\in \partial \Omega$, let $r_{x'}: [0, \tau)\to \bar \Omega$ denote the unit-speed geodesic
starting at $x'$ and normal to $\partial \Omega$. If $x':=\{x_1, \cdots, x_{n-1}\}$  are any local coordinates for
$\partial \Omega$ near $x_0\in \partial \Omega$, we can extend them smoothly to functions on a neighborhood
of $x_0$ in $\Omega$ by letting them be constant along each normal geodesic $r_{x'}$. If we then
define $x_n$ to be the parameter along each $r_{x'}$, it follows easily that $\{x_1, \cdots, x_{n}\}$
form coordinates for $\Omega$ in some neighborhood of $x_0$, which we call the boundary
normal coordinates determined by $\{x_1, \cdots, x_{n-1}\}$. In these coordinates $x_n>0$ in
 $\Omega$, and $\partial \Omega$ is locally characterized by $x_n= 0$. A standard computation shows
that the metric $g$ on $\bar \Omega$ then has the form
 (see p.$\,$1101 of \cite{LU} or p.$\,$532 of \cite{Ta2})
\begin{eqnarray} \label{18/a-1} \; \quad\;\big[g_{jk} (x',x_{n}) \big]_{n\times n} = \begin{bmatrix}
 g_{11} (x',x_{n}) & g_{12} (x',x_{n}) & \cdots & g_{1,n-1} (x',x_{n}) & 0\\
 \cdots\cdots& \cdots\cdots & \cdots &\cdots\cdots  & \cdots\\
 g_{n-1,1} (x',x_{n})  & g_{n-1,2} (x',x_{n}) & \cdots & g_{n-1,n-1} (x',x_{n}) & 0\\
 0& 0& 0& 0&1 \end{bmatrix}.  \end{eqnarray}
       Furthermore, we can take a geodesic normal coordinate system for $(\partial \Omega, g|_{\partial \Omega})$ centered at $x_0=0$, with respect to $e_1, \cdots, e_{n-1}$, where  $e_1, \cdots, e_{n-1}$ are the principal curvature vectors. As Riemann showed, one has (see p.$\,$555 of \cite{Ta2}, or \cite{Spi2})
            \begin{eqnarray} \label{18/7/14/1} \begin{split}& g_{jk}(x_0)= \delta_{jk}, \; \; \frac{\partial g_{jk}}{\partial x_l}(x_0)
 =0  \;\;  \mbox{for all} \;\; 1\le j,k,l \le n-1,\\
 & \frac{1}{2}\,\frac{\partial g_{jk}}{\partial x_n} (x_0) =\kappa_k\delta_{jk}  \;\;  \mbox{for all} \;\; 1\le j,k \le n-1,\end{split}\quad \qquad \qquad \quad
 \end{eqnarray}
where  $\kappa_1,\cdots, \kappa_{n-1}$ are the principal curvatures of $\partial \Omega$ at point $x_0=0$.
Under this normal coordinates, we take $-\boldsymbol{\nu}(x)=[0,\cdots, 0,1]^t$.
By (\ref{18/a-1}) we immediately see that the inverse of metric tensor $g$ in the boundary normal coordinates has form:
    \begin{eqnarray*} g^{-1}(x',x_n) =\begin{bmatrix} g^{11}(x', x_n) & \cdots & g^{1,n-1} (x', x_n)& 0 \\
    \cdots & \cdots  & \cdots  & \cdots\\
    g^{n-1, 1}(x',x_n) & \cdots & g^{n-1,n-1}(x',x_n)& 0\\
    0&\cdots \cdots &0 &1\end{bmatrix}. \end{eqnarray*}
 Note that under  the boundary normal coordinates,  we have $\Gamma_{nk}^n=0$ and $\Gamma_{nn}^l  =0$. Thus, in the boundary normal coordinates, the above system (\ref{2022.6.27-10}) of equations can be written as (see \cite{Liu5})
 \begin{eqnarray} \label{200425-1}
 \frac{\partial^2}{\partial x_n^2} \mathbf{I}_{n+1} +\mathbf{B} \frac{\partial }{\partial x_n} +\mathbf{C}=0, \end{eqnarray}
 where
  \begin{align} & \mathbf{B}:=\begin{bmatrix}
\begin{BMAT}(@, 15pt, 15pt){c.c}{c.c}
  \! \left[\Gamma_{n\beta}^\beta \delta_{jk}    \right]_{n\times n}   &  \left[\!4 \mu^{-1} (\mu+\rho)^{\frac{1}{2}} \Gamma_{\beta n}^j  g^{\alpha\beta}\frac{\partial }{\partial x_\alpha} \!\right]_{n\!\times \!1} \\
          \left[0 \right]_{1\!\times n}  & \Gamma_{n\beta}^\beta  \end{BMAT}
\end{bmatrix}\\
&       +\begin{bmatrix}
\begin{BMAT}(@, 15pt, 15pt){c.c}{c.c}
  \! \left[\delta_{kn} \mu (\mu^{-1})^{;j}\!+\! \frac{\delta_{jk} \rho}{\mu(\mu\!+\!\rho)}\frac{\partial \mu}{\partial x_n}
  \!+\! 2\Gamma_{kn}^j
   \right]_{\!n\times \!n}   &\left[\!\frac{(\mu\!+\!\rho)^{\frac{1}{2}}}{\mu}\!\!\left\{\!\!\!\begin{array} {ll} 2 R^j_n  - 2\mu  \frac{\partial((\mu^{-1})^{;j}\!)}{\partial x_n} - 2\mu \Gamma^j_{\alpha n} (\mu^{-1})^{;\alpha} \\
\!+  ( g^{\alpha\beta} \frac{\partial \Gamma_{\!n\alpha}^j}{\partial x_\beta} \!-\! g^{\alpha\gamma} \Gamma_{\!\beta\gamma}^j \Gamma_{\!n\alpha}^\beta \! \!-\! g^{\alpha\beta} \Gamma_{\!n\gamma}^j \Gamma_{\!\alpha\beta}^\gamma ) \\
 - 2 g^{\alpha \gamma} g^{\beta \sigma} \Gamma_{\gamma \sigma}^j \Gamma_{\alpha\beta}^n
  \end{array}\!\!\!\! \right\}\!\right]_{\!n\!\times \!1}
        \\
     \left[\delta_{nk} \,\mu(\mu+\rho)^{-\frac{1}{2}}\right]_{1\!\times \!n}  &    0  \end{BMAT}
\end{bmatrix}\nonumber \end{align}
 and \begin{eqnarray}\label{2022.6.26-8} \mathbf{C}=\mathbf{C}_2+\mathbf{C}_1+\mathbf{C}_0,\end{eqnarray}
 where
 \begin{align} \label{ 200425-3}& \mathbf{C}_2=g^{\alpha\beta} \frac{\partial^2}{\partial x_\alpha \partial x_\beta} \mathbf{I}_{n+1} +\begin{bmatrix}
\begin{BMAT}(@, 15pt, 15pt){c.c}{c.c}
   \!\left[0\right]_{n\times n}   &   \! \left[2 \mu^{-1} (\mu+\rho)^{\frac{1}{2}} \Gamma_{\gamma \sigma}^j g^{\alpha \gamma} g^{\beta \sigma} \frac{\partial^2}{\partial x_\alpha\partial x_\beta} \right]_{n\times 1}\\
   \left[0 \right]_{1\times n}&  0  \end{BMAT}
\end{bmatrix},  \end{align}
\begin{align}\label{200425-4}  & \,\mathbf{C}_1:=
  \Big( \big( g^{\alpha\beta} \Gamma_{\alpha\gamma}^\gamma +\frac{\partial g^{\alpha\beta}}{\partial x_\alpha} \big)\frac{\partial }{\partial x_\beta} \Big) \mathbf{I}_{n+1}\qquad \qquad \qquad \qquad\qquad \qquad \qquad \qquad\;\;\;\; \end{align}
\begin{align*} &       \!+\!\begin{small}\begin{bmatrix}
\begin{BMAT}(@, 5pt, 5pt){c.c}{c.c}
  \! \!\left[\!\begin{array}{ll}(\!1\!-\!\delta_{nk}\!)\mu(\mu^{\!-\!1}\!)^{;j} \!\frac{\partial }{\partial x_k} \\
  +
  \big(\!\frac{\delta_{jk} \rho}{\mu(\mu\!+\!\rho)}\frac{\partial \mu}{\partial x_\alpha}
  \!+\! 2  \Gamma^j_{\!k\alpha}\!\big)g^{\alpha\beta} \frac{\partial }{\partial x_\beta} \end{array} \! \right]_{\!n\!\times \!n}   \! &  \! \left[\!\frac{(\mu\!+\!\rho)^{\frac{1}{2}}}{\mu}\!\!\left\{\!\!\!\begin{array} {ll} \!2 R^j_\alpha g^{\alpha\beta}  \!\!-\! 2\mu g^{\alpha\beta}  \frac{\partial((\mu^{\!-1}\!)^{;j}\!)}{\partial x_\alpha}\\
  \! -\! 2\mu g^{\alpha\beta} \Gamma^j_{s\alpha } \!(\mu^{\!-1})^{;s}  \!
 -\! 2\Gamma_{sh}^j   g^{sr}\! g^{hm} \Gamma_{rm}^\beta\\
  \!-\!( g^{mr}\! \frac{\partial \Gamma_{\!\alpha r}^j}{\partial x_m}
  \!-\!\! g^{mr} \Gamma_{\!hr}^j \Gamma_{\!\alpha m}^h \!
 \!-\!\! g^{mr} \!\Gamma_{\!\alpha h}^j\! \Gamma_{\!mr}^h )g^{\alpha\!\beta} \\
 \end{array}\!\!\!\!\! \right\}\!\!\frac{\partial }{\partial x_\beta} \!\! \right]_{\!n\!\times \!1}\!\\
     \left[(1-\delta_{nk})\mu(\mu+\rho)^{-\frac{1}{2}}  \frac{\partial }{\partial x_k}\right]_{1\times n} & 0\!\end{BMAT}\!
\end{bmatrix}\end{small},\nonumber
  \end{align*}
 \begin{align} \label{200425-6} & \mathbf{C}_0:=
 \begin{bmatrix}
\begin{BMAT}(@, 15pt, 15pt){c.c}{c.c}
   \!\left[\!\Big((\mu+\rho)^{\!\frac{1}{2}}\Delta_g ((\mu+\rho)^{\!-\frac{1}{2}}) +\mu^{\!-\!1} (\mu\!+\!\rho)^{\frac{1}{2}}   \frac{\partial \mu}{\partial x_l}   ((\mu+\rho)^{\!-\frac{1}{2}})^{;l}\Big) \delta_{jk}  \right]_{n\times n}   &   \! \left[\, 0\,\right]_{n\times 1}\\
   \left[0 \right]_{1\times n}&  -\mu (\Delta_g \mu^{-1})
  \end{BMAT}
\end{bmatrix} \qquad \end{align}
\begin{align*}
\!\!&          \!  +\!\begin{bmatrix}
\begin{BMAT}(@, 1pt, 1pt){c.c}{c.c}
  \! \!\left[\!\!\!\left. \begin{array}{ll}\! \mu(\mu^{\!-\!1}\!)^{;j} \Gamma_{\!kl}^l \!+\!\mu(\mu\!+\!\rho)^{\!\frac{1}{2}}(\mu^{\!-\!1})^{;j} \frac{\partial ((\mu\!+\!\rho)^{\!-\!\frac{1}{2}}\!)}{\partial x_k}\\
 \! \!+\! \big( g^{ml} \frac{\partial \Gamma_{\!kl}^j}{\partial x_m}\! + \!g^{ml}\Gamma_{\!hl}^j \Gamma_{\!km}^h \!\!-\! g^{ml} \Gamma_{\!kh}^j \Gamma_{\!ml}^h \big)\\
   \frac{\rho}{\mu(\mu\!+\!\rho)} g^{ml} \Gamma_{\!km}^j \frac{\partial \mu}{\partial x_l} \!+\!R^j_k \!-\! \mu^{\!
-\!1}\big( \!\frac{\partial \mu}{\partial x_k}\big)^{;j}\end{array} \!\!
\right.\!\! \! \right]_{\!n\!\times \!n}   &   \! \left[\!(\mu\!+\!\rho)^{\!\frac{1}{2}}\!\!\left\{\!\!\! \begin{array}{ll} \!\!-2 ( \Delta_g \mu^{\!-\!1})^{;j} - 2  R_l^j \,(  \mu^{\!-1})^{;l}  \\ \!\!-2  g^{ml} \Gamma_{\!sm}^j \frac{\partial((\mu^{\!-\!1})^{;s}\!)}{\partial x_l}
  -2 \mu^{\!-\!1}\frac{\partial \mu}{\partial x_m} (\mu^{\!-\!1}\!)^{;m;j}\\
    \!\!- \!\big(\! g^{ml} \frac{\partial \Gamma^j_{\!sl}}{\partial x_m}\! + \!g^{ml} \Gamma^j_{\!hl}\Gamma^h_{\!sm} \!\!-\! g^{ml} \Gamma_{\!sh}^j\Gamma^h_{\!ml}\! \big) \!(u^{\!-\!1}\!)^{;s}\end{array}\!\!\!\!\!\right\}\!\!\right]_{\!n\!\times \!1}\\
   \left[\mu(\mu+\rho)^{-\frac{1}{2}} \Gamma_{kl}^l + \mu \frac{\partial ((\mu+\rho)^{-\frac{1}{2}}\!)}{\partial x_k} \! \right]_{\!1\times n}&  0
  \end{BMAT}
\!\end{bmatrix}.\nonumber
 \end{align*}

\noindent Throughout this paper, we denote $\sqrt{-1}={i}$. By applying the method of factorization, we can get a pseudodifferential operator $Q(x, D_{x'})$ of order one in $x'$ depending smoothly on $x_n$ such that \begin{eqnarray*} \label{19.3.19-1:}  \frac{\partial^2}{\partial x_n^2} \mathbf{I}_{n+1} +\mathbf{B}\frac{\partial }{\partial x_n} +\mathbf{C} = \left( \frac{\partial}{\partial x_n} \,\mathbf{I}_{n+1} + \mathbf{B} -\mathbf{Q} \right) \left(\frac{\partial }{\partial x_n}\,\mathbf{I}_{n+1}  +\mathbf{Q}\right), \end{eqnarray*}
 modulo a smoothing operator, where $D_{x'}=(D_{x_1}, \cdots, D_{x_{n-1}})$, $\,D_{x_j}=\frac{1}{{i}}\,\frac{\partial }{\partial x_j}$.
 Let $\mathbf{q}(x, \xi')$, $\mathbf{b}(x,\xi')$ and $\mathbf{c}(x,\xi')$ be the full symbols of $\mathbf{Q}$ and $\mathbf{B}$ and $\mathbf{C}$, respectively. Clearly,  $\mathbf{q}(x, \xi') \sim \sum_{j=0}^\infty \mathbf{q}_{1-j} (x, \xi')$,   $\;\mathbf{b}(x,\xi')=\mathbf{b}_1(x,\xi')+\mathbf{b}_0(x,\xi')$ and $\mathbf{c}(x, \xi')= \mathbf{c}_2(x,\xi') +\mathbf{c}_1(x, \xi') +\mathbf{c}_0(x, \xi')$,
    where \begin{eqnarray}
  \label{19.3.19-4'} & \mathbf{b}_1(x, \xi') = \begin{bmatrix}
\begin{BMAT}(@, 15pt, 15pt){c.c}{c.c}
  \! \left[0   \right]_{n\times n}   &  \left[4{i} \mu^{-1} (\mu+\rho)^{\frac{1}{2}}\Gamma_{\beta n}^j g^{\alpha\beta}  \xi_\alpha\!\right]_{n\!\times \!1} \\
          \left[0 \right]_{1\!\times n}  & 0  \end{BMAT}
\end{bmatrix},\qquad \qquad \qquad \end{eqnarray}
 \begin{eqnarray}\label{19.3.19-5} \!\!\!\!\!\!&& \;\,\;\, \mathbf{b}_0(x, \xi') \!=\! \Gamma_{n\beta}^\beta \mathbf{I}_{n+1} \\
 &&\; +\!\begin{small}\begin{bmatrix}\!
\begin{BMAT}(@, 5pt, 5pt){c.c}{c.c}
  \! \left[\delta_{kn} \mu (\mu^{-1})^{;j}\!+\! \frac{\delta_{jk} \rho}{\mu(\mu\!+\!\rho)} \frac{\partial \mu}{\partial x_n}\!
   +\! 2\Gamma_{\!kn}^j \!
  \right]_{\!n\!\times\! n}   &\left[\!\frac{(\mu\!+\!\rho)^{\frac{1}{2}}}{\mu}\!\!\left\{\!\!\!\begin{array} {ll} 2 R^j_n  \!-\! 2\mu  \frac{\partial(\mu^{\!-\!1})^{\!;j}}{\partial x_n}\! -\! 2\mu \Gamma^j_{\alpha n} (\mu^{\!-\!1})^{\!;\alpha} \\
\!+ \! ( g^{\alpha\beta} \frac{\partial \Gamma_{\!n\alpha}^j}{\partial x_\beta} \!-\! g^{\alpha\gamma} \Gamma_{\!\beta\gamma}^j \Gamma_{\!n\alpha}^\beta \! \!-\! g^{\alpha\beta} \Gamma_{\!n\gamma}^j \Gamma_{\!\alpha\beta}^\gamma ) \\
 - 2 g^{\alpha \gamma} g^{\beta \sigma} \Gamma_{\gamma \sigma}^j \Gamma_{\!\alpha\beta}^n
  \end{array}\!\!\!\! \right\}\!\right]_{\!n\!\times \!1}
        \\
     \left[\delta_{nk} \,\mu(\mu+\rho)^{-\frac{1}{2}}\right]_{1\!\times \!n}  &    0  \!\end{BMAT}\!
\end{bmatrix}\end{small};\nonumber \end{eqnarray}
  \begin{eqnarray} \label{20200515-5} \mathbf{c}_2(x, \!\xi')= -g^{\alpha\beta} \xi_\alpha \xi_\beta \mathbf{I}_{n+1} +\begin{bmatrix}
\begin{BMAT}(@, 15pt, 15pt){c.c}{c.c}
   \left[0\right]_{n\times n}   &   \! \left[-2 \mu^{-1} (\mu+\rho)^{\frac{1}{2}} \Gamma_{\gamma \sigma}^j g^{\alpha \gamma} g^{\beta \sigma} \xi_\alpha\xi_\beta \right]_{n\times 1}\\
   \left[0 \right]_{1\times n}&  0  \end{BMAT}
\end{bmatrix},\end{eqnarray}
\begin{eqnarray} & \label{20200515-6}  \!\!\!\! \mathbf{c}_1(x, \!\xi')\! =
  {i} \, \big( g^{\alpha\beta} \Gamma_{\alpha\gamma}^\gamma +\frac{\partial g^{\alpha\beta}}{\partial x_\alpha} \big) \xi_\beta \,I_{n+1}\qquad \qquad \qquad \qquad \qquad \qquad \qquad \;\;\quad\quad \quad \;\, \;\; \end{eqnarray}
  \begin{align*}&  \!+\! \begin{small}\begin{bmatrix}
\begin{BMAT}(@, 5pt, 5pt){c.c}{c.c}
  \! \!\left[\!\begin{array}{ll} {i}(\!1\!-\!\delta_{nk}\!)\mu(\mu^{\!-\!1}\!)^{;j} \!\xi_k \\
  +  {i}\big(\!\frac{\delta_{jk} \rho}{\mu(\mu\!+\!\rho)}\frac{\partial \mu}{\partial x_\alpha}
  + 2  \Gamma^j_{\!k\alpha}\!\big)g^{\alpha\beta} \xi_\beta\end{array}  \! \right]_{\!n\!\times \!n}   \! &  \! \left[\!\frac{(\mu\!+\!\rho)^{\frac{1}{2}} {i}}{\mu}\!\left\{\!\!\!\begin{array} {ll} \!2 R^j_\alpha g^{\alpha\beta}  \!\!-\! 2\mu g^{\alpha\beta}  \frac{\partial(\mu^{\!-1}\!)^{;j}}{\partial x_\alpha}\\
  \! -\! 2\mu g^{\alpha\beta} \Gamma^j_{s\alpha } \!(\mu^{\!-1})^{;s}  \!
 -\! 2\Gamma_{sh}^j   g^{sr}\! g^{hm} \Gamma_{rm}^\beta\\
  \!-\!( g^{mr}\! \frac{\partial \Gamma_{\!\alpha r}^j}{\partial x_m}
  \!-\!\! g^{mr} \Gamma_{\!hr}^j \Gamma_{\!\alpha m}^h \!
 \!-\!\! g^{mr} \!\Gamma_{\!\alpha h}^j\! \Gamma_{\!mr}^h )g^{\alpha\!\beta} \\
 \end{array}\!\!\!\!\! \right\}\!\xi_\beta \!\! \right]_{\!n\!\times \!1}\!\\
     \left[{i} (1-\delta_{nk})\mu(\mu+\rho)^{-\frac{1}{2}}  \xi_k\right]_{1\times n} & 0\!\end{BMAT}\!
\end{bmatrix}\end{small},\nonumber\end{align*}
 \begin{eqnarray} && \quad\,\;\mathbf{c}_0(x, \xi') \!=
 \begin{small}\begin{bmatrix}
\begin{BMAT}(@, 15pt, 15pt){c.c}{c.c}
   \!\left[\!\Big((\mu+\rho)^{\!\frac{1}{2}}\Delta_g ((\mu+\rho)^{\!-\frac{1}{2}}) +\mu^{\!-\!1} (\mu\!+\!\rho)^{\frac{1}{2}}   \frac{\partial \mu}{\partial x_l}   ((\mu+\rho)^{\!-\frac{1}{2}})^{;l}\Big) \delta_{jk}  \right]_{n\times n}   &   \! \left[\, 0\,\right]_{n\times 1}\\
   \left[0 \right]_{1\times n}&  -\mu (\Delta_g \mu^{-1})
  \end{BMAT}
\end{bmatrix}\end{small} \qquad \end{eqnarray}
\begin{align*}
\!\!&          \!  +\!\begin{small}\begin{bmatrix}
\begin{BMAT}(@, 1pt, 1pt){c.c}{c.c}
  \! \!\left[\!\!\!\left. \begin{array}{ll}\! \mu(\mu^{\!-\!1}\!)^{;j} \Gamma_{\!kl}^l \!+\!\mu(\mu\!+\!\rho)^{\!\frac{1}{2}}(\mu^{\!-\!1})^{;j} \frac{\partial (\mu\!+\!\rho)^{\!-\!\frac{1}{2}}}{\partial x_k}\\
 \! \!+\! \big( g^{ml} \frac{\partial \Gamma_{\!kl}^j}{\partial x_m}\! + \!g^{ml}\Gamma_{\!hl}^j \Gamma_{\!km}^h \!\!-\! g^{ml} \Gamma_{\!kh}^j \Gamma_{\!ml}^h \big)\\
   \frac{\rho}{\mu(\mu\!+\!\rho)} g^{ml} \Gamma_{\!km}^j \frac{\partial \mu}{\partial x_l} \!+\!R^j_k \!-\! \mu^{\!
-\!1}\big( \!\frac{\partial \mu}{\partial x_k}\big)^{;j}\end{array} \!\!
\right.\!\! \! \right]_{\!n\!\times \!n}   &   \! \left[\!(\mu\!+\!\rho)^{\!\frac{1}{2}}\!\!\left\{\!\!\! \begin{array}{ll} \!\!-2 ( \Delta_g \mu^{\!-\!1})^{;j} - 2  R_l^j \,(  \mu^{\!-1})^{;l}  \\ \!\!-2  g^{ml} \Gamma_{\!sm}^j \frac{\partial((\mu^{\!-\!1})^{;s})}{\partial x_l}
  -2 \mu^{\!-\!1}\frac{\partial \mu}{\partial x_m} (\mu^{\!-\!1}\!)^{;m;j}\\
    \!\!- \!\big(\! g^{ml} \frac{\partial \Gamma^j_{\!sl}}{\partial x_m}\! + \!g^{ml} \Gamma^j_{\!hl}\Gamma^h_{\!sm} \!\!-\! g^{ml} \Gamma_{\!sh}^j\Gamma^h_{\!ml}\! \big) \!(u^{\!-\!1}\!)^{;s}\end{array}\!\!\!\!\!\right\}\!\!\right]_{\!n\!\times \!1}\\
   \left[\mu(\mu+\rho)^{-\frac{1}{2}} \Gamma_{kl}^l + \mu \frac{\partial (\mu+\rho)^{-\frac{1}{2}}}{\partial x_k} \! \right]_{\!1\times n}&  0
  \end{BMAT}
\!\end{bmatrix}\end{small}.\nonumber\end{align*}

In \cite{Liu5}, the  full symbol $\mathbf{q}\sim \sum_{l=0}^\infty \mathbf{q}_{1-l}$ of $\mathbf{Q}$ has been obtained:
 \begin{eqnarray}\label{2022.6.23-4}  \boldsymbol{q}_1\! =\sqrt{ g^{\alpha\beta} \xi_\alpha\xi_\beta}\,  \mathbf{I}_{n+1}\!+\!  \mu^{\!-\!1}\! (\mu\!+\!\rho)^{\!\frac{1}{2}} \! \begin{bmatrix}
\begin{BMAT}(@, 15pt, 15pt){c.c}{c.c}
   \!\left[0\right]_{n\times n}   &   \! \left[ 2{i}\, \Gamma_{\beta n}^j  g^{\beta \alpha}\xi_\alpha+
   \frac{\Gamma_{\gamma \sigma}^j g^{\alpha \gamma} g^{\beta \sigma} \xi_\alpha\xi_\beta}{\sqrt{ g^{\alpha\beta}\xi_\alpha\xi_\beta }}\right]_{\!n\times 1}\\
   \left[0 \right]_{1\times n}&  0  \end{BMAT}
\end{bmatrix},\end{eqnarray}
which is a positive-definite matrix.
\begin{eqnarray}&&\label{2022.6.23-5}   \\
&&\mathbf{q}_0=  \frac{1}{2  \sqrt{ g^{\alpha\beta} \xi_\alpha \xi_\beta}} \mathbf{E}_1 -\frac{(\mu+\rho)^{\!\frac{1}{2}}}{ 4\mu g^{\alpha\beta} \xi_\alpha\xi_\beta}  \mathbf{A}_1 \mathbf{E}_1 -\frac{(\mu+\rho)^{\!\frac{1}{2}}}{ 4\mu g^{\alpha\beta} \xi_\alpha\xi_\beta} \mathbf{E}_1 \mathbf{A}_2 +    \frac{(\mu+\rho)}{4\mu^2 \big(g^{\alpha\beta} \xi_\alpha\xi_\beta\big)^{\frac{3}{2}} }  \! \mathbf{A}_1 \mathbf{E}_1 \mathbf{A}_2,\nonumber\end{eqnarray}
where  \begin{eqnarray} \label{19.3.24-1} \mathbf{A}_1 &=&  \begin{small} \begin{bmatrix}
\begin{BMAT}(@, 15pt, 15pt){c.c}{c.c}
   \!\left[0\right]_{n\times n}   &   \! \left[ -2{i}\,  \Gamma_{\beta n}^jg^{\beta \alpha} \xi_\alpha\!+\!
   \frac{1}{\sqrt{ g^{\alpha\beta}\xi_\alpha\xi_\beta }}\Gamma_{\gamma \sigma}^j g^{\alpha \gamma} g^{\beta \sigma} \xi_\alpha\xi_\beta\right]_{\!n\times 1}\\
   \left[0 \right]_{1\times n}&  0  \end{BMAT}
\end{bmatrix} \end{small} \\
    \mathbf{A}_2^t &=&\begin{small}\begin{bmatrix}
\begin{BMAT}(@, 15pt, 15pt){c.c}{c.c}
   \!\left[0\right]_{n\times n}   &   \! \left[0 \right]_{n\times 1}\\
   \left[2{i}\,  \Gamma_{\beta n}^k g^{\beta \alpha}\xi_\alpha\!+\!
   \frac{1}{\sqrt{ g^{\alpha\beta}\xi_\alpha\xi_\beta }}\Gamma_{\gamma \sigma}^k g^{\alpha \gamma} g^{\beta \sigma} \xi_\alpha\xi_\beta \right]_{1\times n}&  0  \end{BMAT}
\end{bmatrix},\end{small}  \end{eqnarray}
and \begin{eqnarray} \label{19.3.23-1} \mathbf{E}_1:= i\sum_{l=1}^{n-1} \frac{\partial \mathbf{q}_1}{\partial \xi_l} \, \frac{\partial \mathbf{q}_1}{\partial x_l} + \mathbf{b}_0\mathbf{q}_1 -i \sum_{l=1}^{n-1} \frac{\partial \mathbf{b}_1}{\partial \xi_l}\, \frac{\partial \mathbf{q}_1}{\partial x_l} +\frac{\partial \mathbf{q}_1}{\partial x_n}  -\mathbf{c}_1,\end{eqnarray}
 and $\mathbf{b}_0$ and $\mathbf{c}_1$ are given in (\ref{19.3.19-5}) and (\ref{20200515-6}).

 Furthermore,
 \begin{eqnarray} \label{19.6.1-1}&& \mathbf{E}_0:=-\mathbf{q}_0^2 +i\sum_{l=1}^{n-1} \big(\frac{\partial \mathbf{q}_1}{\partial \xi_l} \frac{\partial \mathbf{q}_0}{\partial x_l} +\frac{\partial \mathbf{q}_0}{\partial \xi_l}\, \frac{\partial \mathbf{q}_1}{\partial x_l}\big) \\
 && \qquad \quad + \frac{1}{2} \sum_{l,\gamma=1}^{n-1} \frac{\partial^2 \mathbf{q}_1}{\partial \xi_l\partial \xi_\gamma}   \, \frac{\partial^2 \mathbf{q}_1}{\partial x_l \partial x_\gamma}   +\mathbf{b}_0 \mathbf{q}_0-i\sum_{l=1}^{n-1} \frac{\partial \mathbf{b}_1}{\partial \xi_l} \, \frac{\partial \mathbf{q}_0}{\partial x_l} +\frac{\partial \mathbf{q}_0}{\partial x_n} -\mathbf{c}_0.\nonumber\end{eqnarray}
Generally, for $m\ge 1$ we get
 \begin{eqnarray} \label{19.6.1-2} \mathbf{E}_{-m}:= \begin{small}\sum\limits_{\underset{ |\vartheta| = j+k+m}{-m\le j,k\le 1}} \frac{(-i)^{|\vartheta|}}{\vartheta!}\end{small} ( \partial_{\xi'}^{\vartheta} \mathbf{q}_j ) ( \partial_{x'}^{\vartheta} \mathbf{q}_k) +\mathbf{b}_0 \mathbf{q}_{-m}  -i\sum_{l=1}^{n-1} \frac{\partial \mathbf{b}_1}{\partial \xi_l}  \frac{\partial \mathbf{q}_{-m}}{\partial x_l} +\frac{\partial \mathbf{q}_{-m}}{\partial x_n}. \end{eqnarray}
Put \begin{eqnarray} \label{200427-6,} && \mathbf{X}= \frac{1}{2  \sqrt{ g^{\alpha\beta} \xi_\alpha \xi_\beta}} \mathbf{E}  -\frac{(\mu+\rho)^{\!\frac{1}{2}}}{ 4\mu g^{\alpha\beta} \xi_\alpha\xi_\beta} \mathbf{A}_1 \mathbf{E} \\
&& \qquad \;\, -\frac{(\mu+\rho)^{\!\frac{1}{2}}}{ 4\mu g^{\alpha\beta} \xi_\alpha\xi_\beta} \mathbf{E} \mathbf{A}_2 +  \frac{(\mu+\rho)}{4\mu^2 \big(g^{\alpha\beta} \xi_\alpha\xi_\beta\big)^{\frac{3}{2}} }  \! \mathbf{A}_1 \mathbf{E} \mathbf{A}_2.\nonumber\end{eqnarray}
Replacing the matrices $\mathbf{E}$ and $\mathbf{X}$ by the above $\mathbf{E}_{-m}$ and $\mathbf{q}_{-m-1}$ in (\ref{200427-6,}), respectively, we explicitly get all $\mathbf{q}_{-m-1}$, $m\ge 0$ (see \cite{Liu5}).
It has been proved (see \cite{Liu5}) that
 in the local boundary normal coordinates, the Dirichlet-to-Neumann map $\boldsymbol{\Xi}$ associated with (\ref {20200524-4}) can be represented as:
 \begin{eqnarray}\label{19.3.28-10}&&
\frac{\boldsymbol{\partial}}{\boldsymbol{\partial} \mathbf{x_n}}=- \boldsymbol{\Xi}= - \mathbf{Q} \end{eqnarray}
modulo a  smoothing operator.
\vskip 0.26 true cm

\vskip 1.39 true cm

\section{The full symbol of the Dirichlet-to-Neumann map $\boldsymbol{\Lambda}$}

\vskip 0.48 true cm

\noindent{\bf Lemma 4.1.} \ {\it Let $(\Omega,g)$ be a smooth, $n$-dimensional compact Riemannian manifold with smooth boundary. Let $\boldsymbol{\Lambda}$ be the Dirichlet-to-Neumann map associated with Stokes flow. Then the full symbol $\iota(\boldsymbol{\Lambda})$ of the pseudodifferential operator $\boldsymbol{\Lambda}$ has the following asymptotic expansion:
 $$\iota(\boldsymbol{\Lambda})\sim \sum_{m=0}^\infty \boldsymbol{\psi}_{1-m},$$
 where each  $\boldsymbol{\Psi}_{1-m}$ (homogeneous in $\xi'$ of degree $1-m$) can be explicitly calculated (see (\ref{2022.7.1-8}), (\ref{2022.7.5-4}) and (\ref{2022.6.20-6'}) below). In particular,  we have
\begin{eqnarray} \label{2022.7.1-8}&& \boldsymbol{\psi}_{1}= \begin{bmatrix} \psi_1^{11} & \cdots & \psi_1^{1n} \\
\cdots & \cdots & \cdots \\
 \psi_1^{n1} & \cdots & \psi_1^{nn}\end{bmatrix}= 2\mu \sqrt{g^{\eta\theta} \xi_\eta\xi_\theta}\,\;\mathbf{I}_n\\
 &&\label{2022.7.5-4}\boldsymbol{\psi}_{0} =\begin{bmatrix} \psi_0^{11}& \cdots & \psi_0^{1n}\\
\cdots & \cdots & \cdots\\
\psi_0^{n1} &\cdots &\cdots \psi_0^{nn}\!\end{bmatrix} =  \begin{bmatrix} \frac{i\mu \xi_1 \phi_1^1}{g^{\eta\theta} \xi_\eta\xi_\theta} & \cdots & \frac{i\mu \xi_{n-1} \phi_1^1}{g^{\eta\theta} \xi_\eta\xi_\theta} & 0\\
\cdots & \cdots & \cdots\\
\frac{i\mu \xi_1 \phi_1^{n-1}}{g^{\eta\theta} \xi_\eta\xi_\theta} & \cdots & \frac{i\mu \xi_{n-1} \phi_1^{n-1}}{g^{\eta\theta} \xi_\eta\xi_\theta} & 0\\
\frac{i\mu \xi_1 \phi_1^{n}}{g^{\eta\theta} \xi_\eta\xi_\theta} & \cdots & \frac{i\mu \xi_{n-1} \phi_1^{n}}{g^{\eta\theta} \xi_\eta\xi_\theta} & 0
\!\end{bmatrix},  \end{eqnarray}
where \begin{eqnarray*} && \phi_1^{j} = \psi_1^{jj} g^{j\alpha} \frac{\partial \mu^{-1}}{\partial x_\alpha} -\sum\limits_{\gamma=1}^{n-1} \frac{\partial \psi_1^{jj}}{\partial \xi_\gamma} \frac{\partial}{\partial x_\gamma} \big(\mu^{-1} g^{j\alpha}\big) \xi_\alpha -\mu(\mu+\rho)^{-\frac{1}{2}}q_1^{j,n+1} - 2i g^{j\alpha} \xi_\alpha q_0^{n+1,n+1} \\
 && \quad\quad  \;\;- g^{j\alpha} \frac{\partial q_1^{n+1,n+1}}{\partial x_\alpha} + i \frac{\partial g^{j\alpha}}{\partial x_n} \xi_\alpha,\,\;\; 1\le j\le n-1;\\
 &&  \phi_1^n = \psi_1^{nn} \big( \frac{\partial \mu^{-1}}{\partial x_n} +\mu^{-1} q_0^{n+1,n+1} \big) -i \sum\limits_{\gamma=1}^{n-1} \frac{\partial \psi_1^{nn}}{\partial \xi_\gamma} \,\frac{\partial}{\partial x_\gamma} \big(\mu^{-1} q_1^{n+1,n+1}\big)- 2\mu (\mu+\rho)^{-\frac{1}{2}}q_1^{n,n+1} \\
 && \;\; \;\;\quad \;\; + 3\sum_{l=1}^n q_0^{n+1,l} q_1^{l, n+1} -3\frac{\partial q_1^{n+1,n+1}}{\partial x_n} +6 q_1^{n+1,n+1} q_0^{n+1,n+1}\\
&& \quad \quad\;\;  -3i \sum\limits_{\gamma=1}^{n-1} \frac{\partial q_1^{n+1,n+1} }{\partial \xi_\gamma}\, \frac{\partial q_1^{n+1,n+1}}{\partial x_\gamma}
      -\Gamma^{\beta}_{n\beta}  q_1^{n+1,n+1}+ \frac{{i}}{\sqrt{|g|}}\frac{\partial}{\partial x_\alpha}\big(\sqrt{|g|}g^{\alpha\beta}\big)\xi_\beta,\end{eqnarray*}
$q_l^{j,k}$ is the $(j,k)$ entry of $\mathbf{q}_l=[q_l^{jk}]_{(n+1)\times (n+1)}$, ($l=1,0$), and $\mathbf{q}_l$ is given by (\ref{2022.6.23-4}) and (\ref{2022.6.23-5}). }
\vskip 0.68 true cm

\noindent  {\it Proof of Lemma 4.1.} \ It follows from \cite{Liu5} (or \cite{Liu1}) that
\begin{eqnarray} \label{2022.6.9-1}\!\!\!\! &&\;\;\;\;\;\;\;\;\sigma_\mu (\mathbf{u}, p) (-\boldsymbol{\nu}) = 2\mu (\mbox{Def}\, \mathbf{u})^\# (-\boldsymbol{\nu}) -p (-\boldsymbol{\nu})\\ [2.5mm]
\!\!\!\!&&  =\begin{bmatrix}
\begin{BMAT}(@, 15pt, 15pt){c.c}{c.c}
  \! \left[\mu \delta_{jk} \frac{\partial}{\partial x_n} \right]_{(n-1)\times (n-1)}   &  \left[\!\mu g^{j\alpha} \frac{\partial }{\partial x_\alpha}  \!\right]_{(n\!-1)\times \!1} \\
          \left[0 \right]_{1\!\times (n-1)}  & 2\mu \frac{\partial}{\partial x_n}  \end{BMAT}
\end{bmatrix}\!\begin{bmatrix} u^1\\
\cdots \\
u^{n-1} \\
u^n\end{bmatrix} +\begin{bmatrix} 0\\
\cdots \\
0 \\
p\end{bmatrix}\nonumber\\ [1mm]
\!\!\!\!&& =\begin{bmatrix}
\begin{BMAT}(@, 15pt, 15pt){c.c}{c.c}
  \! \left[\mu \delta_{jk} \frac{\partial}{\partial x_n} \right]_{(n-1)\times (n-1)}   &  \left[\!\mu g^{j\alpha} \frac{\partial }{\partial x_\alpha}  \!\right]_{(n\!-1)\times \!1} \\
          \left[0 \right]_{1\!\times (n-1)}  & 2\mu \frac{\partial}{\partial x_n}  \end{BMAT}
\end{bmatrix} \!
\left[
\begin{array}{c;{3pt/3pt}c}
  &  \mu^{-1}g^{1l} \frac{\partial}{\partial x_l} -g^{1l}\frac{\partial \mu^{-1}}{\partial x_l}  \\
   \!\!\! (\mu+\rho)^{-\frac{1}{2}}I_n\!\!   & \cdots\\
   & \mu^{-1}g^{nl} \frac{\partial}{\partial x_l} -g^{nl}\frac{\partial \mu^{-1}}{\partial x_l} \\
 \end{array}\!\!
\right]\!\begin{bmatrix} w^1\\
\cdots \\
w^{n-1} \\
f\end{bmatrix} \nonumber\\ [1.5mm]
\!\!\!\!&&\;\;\;\;+
\left[
\begin{array}{ccc;{3pt/3pt}c}
0 & \cdots & 0 & 0 \\
 \vdots  & \cdots & \vdots  & \vdots \\
0 & \cdots & 0 & 0 \\
(\mu + \rho)^{-\frac{1}{2}}\frac{\partial \mu}{\partial x_{1}} & \cdots & (\mu + \rho)^{-\frac{1}{2}}\frac{\partial \mu}{\partial x_{n}} &   \Delta_{g} + \mu (\Delta_{g} \mu^{-1}) \\
\end{array}
\right]\!\begin{bmatrix} w^1\\
\cdots \\
w^{n-1} \\
f\end{bmatrix} \nonumber
\\ [2.6mm]
\!\!\!\!&& =\!
    \begin{bmatrix}\!\!
        \begin{BMAT}(@, 30pt, 30pt){c.c.c}{c.c}
        H_n I_{n-1} &
        \bigl[g^{j\alpha} H_\alpha\bigr]_{(n-1)\times 1} &
        \biggl[J_j- \mu\frac{\partial }{\partial x_n}\Bigl(g^{j\alpha}\frac{\partial \mu^{-1}}{\partial x_\alpha}\Bigr) - \mu g^{j\alpha}\frac{\partial^2 \mu^{-1}}{\partial x_\alpha x_n}\biggr] \\
        \frac{1}{(\mu+\rho)^{1/2}}\biggl[\frac{\partial \mu}{\partial x_k}\biggr] & 2 H_n \!+\!\frac{1}{(\mu\!+\!\rho)^{\!1/2}}\frac{\partial \mu}{\partial x_n} & 3\frac{\partial^2}{\partial x_n^2} \!+\!\Gamma^\beta_{n\beta} \frac{\partial}{\partial x_n} \!-\! 2\mu \frac{\partial^2 \mu^{-1}}{\partial x_n^2}\! +\! \Delta_{\partial \Omega} \!+\! \mu \Delta_g \mu^{\!-1}
        \end{BMAT}
    \end{bmatrix}\!\begin{bmatrix} w^1\\
\cdots \\
w^{n-1} \\
f\end{bmatrix}, \nonumber
 \end{eqnarray}
where   \begin{align*} & H_j = \mu \Bigl((\mu+\rho)^{-1/2}\frac{\partial }{\partial x_j} + \frac{\partial (\mu+\rho)^{-1/2}}{\partial x_j}\Bigr) \;\, \mbox{for}\;\, 1\le j \le n,\\
  &  J_j=  g^{j\alpha} \big(\frac{\partial^2}{\partial x_n\partial x_\alpha} +\frac{\partial^2}{\partial x_\alpha\partial x_n} \big) + \frac{\partial g^{j\alpha}}{\partial x_n} \frac{\partial }{\partial x_\alpha} \;\, \mbox{for}\;\, 1\le j\le n-1,\\
 &   \Delta_{\partial \Omega}:=\frac{1}{\sqrt{|g|}}\frac{\partial}{\partial x_\alpha}\big(\sqrt{|g|}g^{\alpha\beta} \frac{\partial}{\partial x_\beta}\big).\end{align*}
In the third equality of (\ref{2022.6.9-1}), we have used the transform (\ref{2022.6.25-7}). Note that the matrix-valued pseudodifferential operator $\mathbf{Q}$ can be written as \begin{eqnarray}\label{2022.6.26-9} \mathbf{Q}= \begin{bmatrix} Q^{11} &\cdots& Q^{1,n+1}\\
  \vdots && \vdots \\
  Q^{n+1,1}& \cdots & Q^{n+1,n+1} \end{bmatrix}.\end{eqnarray}
   Replacing $\frac{\partial w^j}{\partial x_n}$ and $\frac{\partial f}{\partial x_n}$ by $-\sum_{l=1}^n Q^{j,l} w^l - Q^{j, n+1} f$ and $-\sum_{l=1}^n Q^{n+1,l} w^l - Q^{n+1, n+1} f$, respectively, we have
 \begin{eqnarray} \label{2022.6.19-1} \sigma (\mathbf{u}, p) (-\boldsymbol{\nu})=\mathbf{M} \begin{bmatrix} w^1\\
\cdots \\
w^{n} \\
f\end{bmatrix},\end{eqnarray}
         where \begin{eqnarray} \label{2022.7.5-1} \mathbf{M}= \begin{bmatrix}
                  M^{11} &\cdots & M^{1,n} & M^{1,n+1}\\
               \cdots& \cdots& \cdots& \cdots \\
               M^{n1} & \cdots & M^{nn}& M^{n,n+1}\end{bmatrix}\end{eqnarray} and
          \begin{align*}&  M^{jk}=-\mu(\mu+\rho)^{-\frac{1}{2}} Q^{jk} -g^{j\alpha}\frac{\partial Q^{n+1,k}}{\partial x_\alpha} -2g^{j\alpha} Q^{n+1,k} \frac{\partial }{\partial x_\alpha} +\delta_{jk} \,\mu\frac{\partial (\mu+\rho)^{-\frac{1}{2}}}{\partial x_n}, \;\; 1\le j, k\le n-1;\\
        & M^{jn} =- \mu(\mu+\rho)^{-\frac{1}{2}} Q^{jn} -g^{j\alpha}\frac{\partial Q^{n+1,n}}{\partial x_\alpha} -2g^{j\alpha} Q^{n+1,n} \frac{\partial }{\partial x_\alpha}\\
          & \qquad \quad + \mu g^{j\alpha} \Big( (\mu+\rho)^{-\frac{1}{2}} \frac{\partial}{\partial x_\alpha} +\frac{\partial (\mu+\rho)^{-\frac{1}{2}}}{\partial x_\alpha} \Big), \;\; 1\le j\le n-1;\\
         & M^{nk}= (\mu+\rho)^{-\frac{1}{2}} \frac{\partial \mu}{\partial x_k}- 2\mu (\mu+\rho)^{-\frac{1}{2}}Q^{n, k} -3 \frac{\partial Q^{n+1, k}}{\partial x_n} +3\sum_{l=1}^n Q^{n+1, l}Q^{l,k} \\
         &\qquad \;\quad  +3 Q^{n+1, n+1} Q^{n+1,k} -\Gamma_{n\beta}^\beta Q^{n+1,k}+2\delta_{nk}\,\mu \frac{\partial (\mu+\rho)^{-\frac{1}{2}}}{\partial x_n}, \;\,\; 1\le k\le n;\\
                        &M^{j,n+1} = - \mu(\mu+\rho)^{-\frac{1}{2}} Q^{j,n+1} -g^{j\alpha}\frac{\partial Q^{n+1,n+1}}{\partial x_\alpha} -2g^{j\alpha} Q^{n+1,n+1} \frac{\partial }{\partial x_\alpha}\\
          & \qquad \quad \quad \,+  \frac{\partial g^{j\alpha} }{\partial x_n} \frac{\partial}{\partial x_\alpha} - \mu \frac{ \partial }{\partial x_n} \big( g^{j\alpha} \frac{\partial \mu^{-1} }{\partial x_\alpha} \big) -\mu g^{j\alpha} \frac{\partial^2\mu^{-1}}{\partial x_\alpha \partial x_n}, \;\; 1\le j\le n-1;   \\
          & M^{n,n+1}=- 2\mu (\mu+\rho)^{-\frac{1}{2}} Q^{n,n+1}  +3\sum_{l=1}^n Q^{n+1, l}Q^{l,n+1} -3 \frac{\partial Q^{n+1,n+1}}{\partial x_n} +3 Q^{n+1, n+1} Q^{n+1,n+1}\\
          &\qquad \quad \;\quad    -\Gamma^\beta_{n\beta} Q^{n+1,n+1}  -2\mu \frac{\partial^2 \mu^{-1}}{\partial x_n^2} +\Delta_{\partial \Omega} +\mu \big(\Delta_g \mu^{-1}\big).     \end{align*}
Since
\begin{eqnarray*} && \begin{bmatrix} u^1\\
\vdots \\
u^n\end{bmatrix} =\left[
\begin{array}{c;{3pt/3pt}c}
  &  \mu^{-1}g^{1l} \frac{\partial}{\partial x_l} -g^{1l}\frac{\partial \mu^{-1}}{\partial x_l}  \\
   \!\!\! (\mu+\rho)^{-\frac{1}{2}}I_n\!\!   & \cdots\\
   & \mu^{-1}g^{nl} \frac{\partial}{\partial x_l} -g^{nl}\frac{\partial \mu^{-1}}{\partial x_l} \\
 \end{array}\!\!
\right]\!\begin{bmatrix} w^1\\
\cdots \\
w^{n} \\
f\end{bmatrix} \nonumber \\ [3mm]
&& =\left[
\begin{array}{c;{3pt/3pt}c}
  &  \mu^{-1}g^{1\alpha} \frac{\partial}{\partial x_\alpha} -g^{1\alpha}\frac{\partial \mu^{-1}}{\partial x_\alpha}  \\  &{} \\
   \!\!\! (\mu+\rho)^{-\frac{1}{2}}I_n\!\!   & \vdots\\
   & \mu^{-1}g^{n-1,\alpha} \frac{\partial}{\partial x_\alpha} -g^{n-1,\alpha}\frac{\partial \mu^{-1}}{\partial x_l} \\
  & \mu^{-1} \frac{\partial}{\partial x_n} -\frac{\partial \mu^{-1}}{\partial x_n}
 \end{array}\!\!
\right]\!\begin{bmatrix} w^1\\
\cdots \\
w^{n} \\
f\end{bmatrix} \nonumber \\ [3mm]
&&=\begin{bmatrix} (\mu+\rho)^{-\frac{1}{2}} w^1 + (\mu^{-1} g^{1\alpha}\frac{\partial}{\partial x_\alpha} -g^{1\alpha}\frac{\partial \mu^{-1}}{\partial x_\alpha})f \\
\cdots  \\
(\mu+\rho)^{-\frac{1}{2}} w^{n-1} + (\mu^{-1} g^{n-1,\alpha}\frac{\partial}{\partial x_\alpha} -g^{n-1,\alpha}\frac{\partial \mu^{-1}}{\partial x_\alpha})f \\
(\mu+\rho)^{-\frac{1}{2}} w^n-  \mu^{-1}( \sum_{l=1}^nQ^{n+1,l}w^l +Q^{n+1,n+1} f)  -\frac{\partial \mu^{-1}}{\partial x_n})f
\end{bmatrix}\\
[3mm]
&& =\mathbf{L} \begin{bmatrix} w^1\\
\cdots \\
w^{n-1} \\
f\end{bmatrix}, \nonumber \end{eqnarray*}
where  \begin{eqnarray*} \mathbf{L}:= \left[\!\!
\begin{array}{ccc;{3pt/3pt}c}
(\mu+\rho)^{-\frac{1}{2}} \delta_{11} & \cdots & (\mu+\rho)^{-\frac{1}{2}} \delta_{1n} & \mu^{-1} g^{1\alpha}\frac{\partial}{\partial x_\alpha} -g^{1\alpha}\frac{\partial \mu^{-1}}{\partial x_\alpha} \\
 \vdots  & \cdots & \vdots  & \vdots \\
(\mu+\rho)^{-\frac{1}{2}} \delta_{n-1,1} & \cdots & (\mu+\rho)^{-\frac{1}{2}} \delta_{n-1,n} & \mu^{-1} g^{n-1,\alpha}\frac{\partial}{\partial x_\alpha} -g^{n-1,\alpha}\frac{\partial \mu^{-1}}{\partial x_\alpha} \\
-\mu^{-1}Q^{n+1, 1}  & \cdots & -\mu^{-1}Q^{n+1,n} +(\mu+\rho)^{-\frac{1}{2}} &   -\mu^{-1} Q^{n+1, n+1} -\frac{\partial \mu^{-1}}{\partial x_n}
\end{array}\!\!\!
\right].
\nonumber \end{eqnarray*}

Our aim is to look for a matrix-valued the pseudodifferential operator \begin{eqnarray*}\boldsymbol{\Psi}:=\begin{bmatrix} \Psi^{11} & \cdots & \Psi^{1n}\\
\cdots & \cdots & \cdots\\
\Psi^{11} & \cdots & \Psi^{1n}\end{bmatrix} \end{eqnarray*}
such that  $\boldsymbol{\Psi} \mathbf{L}=\mathbf{M}$.
Let us write the full symbol of  $\boldsymbol{\Psi}$  as
\begin{eqnarray*} \iota(\boldsymbol{\Psi})\;\sim \;\sum_{m=0}^\infty \boldsymbol{\psi}_{1-m}, \end{eqnarray*}
where \begin{eqnarray*}\boldsymbol{\psi}_{1-m}:= \begin{bmatrix}\psi_{1-m}^{11} & \cdots & \psi_{1-m}^{1n}\\
\cdots & \cdots & \cdots\\
\psi_{1-m}^{11} & \cdots & \psi_{1-m}^{1n}\end{bmatrix}\end{eqnarray*}
is homogeneous  of degree $1-m$ in $\xi$ for $|\xi|>1$.
Combining this, the operator equation $\boldsymbol{\Psi} \mathbf{L}=\mathbf{M}$  and symbol
formula (\ref{2022.7.1-4}) for product of two pseudodifferential operators, we get the full
symbol equation:
\begin{eqnarray}\label{2022.6.10-2} \sum_{\vartheta} \frac{(-i)^{|\vartheta|}}{\vartheta!} (\partial^\vartheta_{\xi'} \boldsymbol{\psi} ) (\partial_{x'}^{\vartheta} \boldsymbol{\it l}) = \sigma(\mathbf{M}),\end{eqnarray}
where $\boldsymbol{\it l}\sim \sum_{m=0}^\infty \boldsymbol{\it l}_{1-m}$ is the full symbol of $\mathbf{L}$, and
\begin{eqnarray*}&&\boldsymbol{\it l}_{1}:= \begin{bmatrix}
        \begin{BMAT}(@, 30pt, 30pt){c.c}{c.c}
       \left[ \,0\,\right]_{(n-1)\times n} &\left[i\mu^{-1}g^{j\alpha}\xi_\alpha\right]_{(n-1)\times 1}\\
        \left[-\mu^{-1}q_1^{n+1,k}\right]_{1\times n}&  -\mu^{-1}q_1^{n+1,n+1}
                \end{BMAT}
    \end{bmatrix},\\
   && \boldsymbol{\it l}_{0}:=\left[\!\!
\begin{array}{cccc;{3pt/3pt}c}
   (\mu+\rho)^{-\frac{1}{2}} & \cdots &0 & 0&  -g^{1\alpha} \frac{\partial \mu^{-1}}{\partial x_\alpha}   \\
  \vdots &  \ddots &\vdots &\vdots & \\
  0 &\cdots & (\mu+\rho)^{-\frac{1}{2}} & 0 & -g^{n-1,\alpha} \frac{\partial \mu^{-1}}{\partial x_\alpha}\\
      - \mu^{\!-\!1}  q_0^{n+1,1} &\cdots & - \mu^{\!-\!1}  q_0^{n+1,n-1} & -\mu^{\!-\!1}  q_0^{n+1,n} \!+ \! (\mu\!+\!\rho)^{-\frac{1}{2}} & - \mu^{\!-\!1}  q_0^{n+1,n+1}\!-\!\frac{\partial \mu^{\!-\!1}}{\partial x_n} \end{array}\!\!
\right],\\
 && \boldsymbol{\it l}_{1-m}:= \left[\!\!
\begin{array}{ccc;{3pt/3pt}c}
 0& \cdots &0 & 0 \\
 \vdots  & \cdots & \vdots  & \vdots \\
0& \cdots & 0 & 0 \\
-\mu^{-1}q_{1-m}^{n+1, 1}  & \cdots & -\mu^{-1}q^{n+1,n}_{1-m} &   -\mu^{-1} q^{n+1, n+1}_{1-m}
\end{array}
\right], \;\;\; m=2, 3, 4, \cdots.
    \end{eqnarray*}
Group the homogeneous terms of degree two in (\ref{2022.6.10-2}) we obtain the matrix equation
\begin{eqnarray*}&& \begin{bmatrix} \psi_1^{11} & \cdots & \psi_1^{1n} \\
\cdots & \cdots & \cdots \\
 \psi_1^{n1} & \cdots & \psi_1^{nn}\end{bmatrix}
 \begin{bmatrix}
        \begin{BMAT}(@, 30pt, 30pt){c.c}{c.c}
       \left[ \,0\,\right]_{(n-1)\times n} &\left[i\mu^{-1}g^{j\alpha}\xi_\alpha\right]_{(n-1)\times 1}\\
        \left[-\mu^{-1}q_1^{n+1,k}\right]_{1\times n}&  -\mu^{-1}q_1^{n+1,n+1}
                \end{BMAT}
    \end{bmatrix}\\
    && =
     \begin{bmatrix}
        \begin{BMAT}(@, 30pt, 30pt){c.c}{c.c}
       \left[ -2i g^{j\alpha} \xi_\alpha q_1^{n+1,k} \right]_{(n-1)\times n} &\left[ -2i g^{j\alpha} \xi_\alpha q_1^{n+1,n+1}\right]_{(n-1)\times 1}\\
        \left[3\sum_{l=1}^n q_1^{n+1,l} q_1^{l,k} + 3q_1^{n+1,n+1} q_1^{n+1,k}\right]_{1\times n} &  3\sum_{l=1}^n q_1^{n+1,l} q_1^{l,n+1} +3 q_1^{n+1,n+1} q_1^{n+1,n+1}-g^{\alpha\beta} \xi_\alpha\xi_\beta
                \end{BMAT}
    \end{bmatrix}
    \end{eqnarray*}
    Noting that $q_1^{n+1,1} =\cdots= q_1^{n+1,n} =0$ and  $q_1^{n+1, n+1}=\sqrt{g^{\alpha \beta \xi_\alpha \xi_\beta}}$, we have
    \begin{eqnarray}\label{2022.6.10-3}&& \begin{bmatrix} \psi_1^{11} & \cdots & \psi_1^{1n} \\
\cdots & \cdots & \cdots \\
 \psi_1^{n1} & \cdots & \psi_1^{nn}\end{bmatrix}
 \begin{bmatrix}
        \begin{BMAT}(@, 30pt, 30pt){c.c}{c.c}
       \left[ \,0\,\right]_{(n-1)\times n} &\left[i\mu^{-1}g^{j\alpha}\xi_\alpha\right]_{(n-1)\times 1}\\
        \left[\,0\,\right]_{1\times n}&  -\mu^{-1}q_1^{n+1,n+1}
                \end{BMAT}
    \end{bmatrix}\\
    && =
     \begin{bmatrix}
        \begin{BMAT}(@, 30pt, 30pt){c.c}{c.c}
       \left[\, 0\,\right]_{(n-1)\times n} &\left[ -2i g^{j\alpha} \xi_\alpha q_1^{n+1,n+1}\right]_{(n-1)\times 1}\\
        \left[\,0\,\right] & 3 q_1^{n+1,n+1} q_1^{n+1,n+1}-g^{\alpha\beta} \xi_\alpha\xi_\beta
                \end{BMAT}
    \end{bmatrix}.\nonumber
    \end{eqnarray}
 It can be verified that (\ref{2022.6.10-3}) has a solution of the following form:
 \begin{eqnarray*} &&\psi_1^{11} =-2\mu q_1^{n+1,n+1}, \,\; \psi_1^{12}=\cdots =\psi_1^{1n}=0,\\[1.5mm]
 &&  \quad\cdots  \cdots \cdots \cdots\\[1.5mm]
  && \psi_1^{n-1,1} = \cdots= \psi_1^{n-1, n-2} =0,\, \; \psi_1^{n-1, n-1}= -2\mu q_1^{n+1,n+1},\, \; \psi_1^{n-1,n}=0,\\
  && \psi_1^{n1}=\cdots =\psi_1^{n,n-1} =0, \,\; \psi_1^{nn} =- \mu \Big(3q_1^{n+1, n+1}-\frac{g^{\alpha\beta} \xi_\alpha\xi_\beta}{q^{n+1,n+1}_1}\Big),\end{eqnarray*}
 i.e.,
  \begin{eqnarray*}\label{2022.6.10-5} \begin{bmatrix} \psi_1^{11} & \cdots & \psi_1^{1n} \\
\cdots & \cdots & \cdots \\
 \psi_1^{n1} & \cdots & \psi_1^{nn}\end{bmatrix}=- 2\mu \sqrt{g^{\alpha\beta}\xi_\alpha\xi_\beta}\,\;\mathbf{I}_n. \end{eqnarray*}
 The terms of degree one in (\ref{2022.6.10-2} ) are
 \begin{align} &\label{2022.7.1-6} \\
  & \begin{bmatrix} \psi_1^{11} & &\Huge{0} & \\
  & \ddots & \\
 \Huge {0}  && \psi_1^{nn}\!\!\!\end{bmatrix}\!
    \left[\!\!
\begin{array}{cccc;{3pt/3pt}c}
   (\mu+\rho)^{-\frac{1}{2}} & \cdots &0 & 0&  -g^{1\alpha} \frac{\partial \mu^{-1}}{\partial x_\alpha}   \\
  \vdots &  \ddots &\vdots &\vdots & \\
  0 &\cdots & (\mu+\rho)^{-\frac{1}{2}} & 0 & -g^{n-1,\alpha} \frac{\partial \mu^{-1}}{\partial x_\alpha}\\
      - \mu^{\!-\!1}  q_0^{n+1,1} &\cdots & - \mu^{\!-\!1}  q_0^{n+1,n-1} & -\mu^{\!-\!1}  q_0^{n+1,n} \!+ \! (\mu\!+\!\rho)^{-\frac{1}{2}} & - \mu^{\!-\!1}  q_0^{n+1,n+1}\!-\!\frac{\partial \mu^{\!-\!1}}{\partial x_n} \end{array}\!\!
\right]\nonumber \\
&\;\;\;\;\; + \begin{bmatrix} \psi_0^{11}& \cdots & \psi_0^{1n}\\
\cdots & \cdots & \cdots\\
\psi_0^{n1} &\cdots &\cdots \psi_0^{nn}\!\end{bmatrix} \left[\!\!
\begin{array}{c;{3pt/3pt}c}
  &   i\mu^{-1}g^{1\alpha} \xi_\alpha   \\
      \left[\,0\,\right]_{n\times n}& \vdots\\
    & i\mu^{-1}g^{n-1,\alpha} \xi_\alpha\\
     & -\mu^{-1} q_1^{n+1,n+1} \end{array}\!\!
\right]\qquad \;\nonumber \\
& \;\;\; \;\; -i \begin{bmatrix} \frac{\partial \psi_1^{11}}{\partial \xi'} & & \Huge{0}& \\
  & \ddots & \\
  \Huge{0}& & \frac{\partial \psi_1^{nn}}{\partial \xi'}\!\!\end{bmatrix}\left[\!\!
\begin{array}{c;{3pt/3pt}c}
  &   i\frac{\partial }{\partial x'}(\mu^{-1}g^{1\alpha}) \xi_\alpha   \\
     \left[\,0\,\right]_{n\times n}& \vdots\\
    &  i\frac{\partial }{\partial x'}(\mu^{-1}g^{n-1,\alpha}) \xi_\alpha\\
     & -\frac{\partial }{\partial x'}(\mu^{-1}q_1^{n+1,n+1}) \end{array}\!\!
\right]\nonumber \\
&  =\left[\!\!
\begin{array}{cccc;{3pt/3pt}c}
   m_1^{11} & \cdots & m_1^{1,n-1} & m_1^{1n} & m_1^{1,n+1}    \\
     \cdots & \cdots & \cdots & \cdots& \cdots \\
     m_1^{n-1,1} & \cdots & m_1^{n-1,n-1} & m_1^{n-1,n} & m_1^{n-1,n+1}  \\
     m_1^{n1} & \cdots & m_1^{n,n-1} & m_1^{nn} & m_1^{n,n+1}   \end{array}\!\!
\right],\quad \quad\nonumber
  \end{align}
  where \begin{eqnarray*}
  && m_1^{jk}=-\mu (\mu+\rho)^{-\frac{1}{2}} q_1^{jk} -2i g^{j\alpha} \xi_\alpha q_0^{n+1,k} - g^{j\alpha} \frac{\partial q_1^{n+1,k}}{\partial x_\alpha}, \;\; \,  1\le j,k\le n-1;\\
   &&   m_1^{jn}=-\mu (\mu+\rho)^{-\frac{1}{2}} q_1^{jn} -2i g^{j\alpha} \xi_\alpha q_0^{n+1,n} - g^{j\alpha} \frac{\partial q_1^{n+1,n}}{\partial x_\alpha}+i\mu(\mu+\rho)^{-\frac{1}{2}}g^{j\alpha}\xi_\alpha, \;\; \,  1\le j\le n-1;\\
    &&  m_1^{nk} =- 2\mu (\mu+\rho)^{-\frac{1}{2}} q^{n,k}_1 -3\frac{\partial q_1^{n+1,k}}{\partial x_n} +3\sum_{l=1}^n \big(q_0^{n+1, l}q_1^{l,k} +q_1^{n+1,l}q_0^{l,k} -i \sum\limits_{\gamma=1}^{n-1}\frac{\partial q_1^{n+1,l}}{\partial \xi_\gamma} \frac{\partial q_1^{l,k}}{\partial x_\gamma} \big) \\
    && \quad \quad  + 3q_0^{n+1, n+1}q_1^{n+1,k} +3q_1^{n+1,n+1}q_0^{n+1,k} -3i \sum\limits_{\gamma=1}^{n-1}\frac{\partial q_1^{n+1,n+1}}{\partial \xi_\gamma} \, \frac{\partial q_1^{n+1,k}}{\partial x_\gamma}- \Gamma_{n\beta}^\beta q_1^{n+1, k}, \;\; 1\le k\le n;\\
     &&  m_1^{j,n+1} = -\mu (\mu+\rho)^{-\frac{1}{2}} q_1^{j, n+1} -2i g^{j\alpha}\xi_\alpha q_0^{n+1, n+1} -  g^{j\alpha} \frac{\partial q_1^{n+1, n+1}}{\partial x_\alpha} +i \frac{\partial g^{j \alpha}}{\partial x_n} \xi_\alpha, \;\, 1\le j\le n-1;\\
     && m_1^{n,n+1} =- 2\mu (\mu+\rho)^{-\frac{1}{2}} q_1^{n, n+1} \!+ \!3\sum_{l=1}^n \Big( q_1^{n+1,l} q_0^{l, n+1} \!+\!q_0^{n+1,l} q_1^{l, n+1} \!-i\sum\limits_{\gamma=1}^{n-1}\frac{\partial q_1^{n+1,l}}{\partial \xi_\gamma}\,\frac{\partial q_1^{l,n+1}}{\partial x_\gamma}\Big)  \\
      &&  \qquad -3\frac{\partial q_1^{n+1,n+1}}{\partial x_n} \!+\! 3q_1^{n+1,n+1} q_0^{n+1,n+1}\! +\!3 q_0^{n+1,n+1}q_1^{n+1,n+1}\! -\!3 i \sum\limits_{\gamma=1}^{n-1}\frac{\partial q_1^{n+1,n+1}}{\partial \xi_\gamma} \frac{\partial q_1^{n+1,n+1}}{\partial x_\gamma} \\
      &&\qquad - \Gamma_{n\beta}^\beta q_1^{n+1, n+1} + \frac{{i}}{\sqrt{|g|}} \frac{\partial}{\partial x_\alpha} \big( \sqrt{|g|}\, g^{\alpha\beta}\big)\xi_\beta,
        \end{eqnarray*}
i.e.,
\begin{eqnarray*} && \begin{bmatrix} \psi_0^{11}& \cdots & \psi_0^{1n}\\
\cdots & \cdots & \cdots\\
\psi_0^{n1} &\cdots &\cdots \psi_0^{nn}\!\end{bmatrix} \left[\!\!
\begin{array}{c;{3pt/3pt}c}
  &   {i}\mu^{-1}g^{1\alpha} \xi_\alpha   \\
     \left[\,0\,\right]_{n\times n}& \vdots\\
    &{i}\mu^{-1}g^{n-1,\alpha} \xi_\alpha\\
     & -\mu^{-1} q_1^{n+1,n+1} \end{array}\!\!
\right]\\ [2.5mm]
&&\!\!\! =
\begin{small} \left[\!\!
\begin{array}{cccc;{3pt/3pt}c}
\!-  \psi_1^{11} (\mu\!+\!\rho)^{\!-\frac{1}{2}} & \cdots &0 & 0& \! \psi_1^{11} g^{1\alpha} \frac{\partial \mu^{-1}}{\partial x_\alpha}   \\
  \!\vdots &  \ddots &0 &0 & \\
  \!0 &\cdots & -\psi_1^{n\!-\!1,n\!-\!1}(\mu\!+\!\rho)^{-\frac{1}{2}} & 0 & \psi_1^{n\!-\!1,n\!-\!1} g^{n-1,\alpha} \frac{\partial \mu^{-1}}{\partial x_\alpha}\\
    \! \psi_1^{nn} \mu^{\!-\!1}  q_0^{n+1,1} &\cdots & \psi_1^{nn} \mu^{\!-\!1}  q_0^{n+1,n-1} & - \psi_1^{nn} \big(-\mu^{\!-\!1}  q_0^{n\!+\!1,n} \!+ \! (\mu\!+\!\rho)^{-\frac{1}{2}}\big) & \! \psi_1^{nn} (\frac{\partial \mu^{\!-\!1}}{\partial x_n}\!+\! \mu^{\!-\!1}  q_0^{n+1,n+1}) \end{array}\!\!
\right]\end{small}
\\ [2.5mm]
&& \; \;\; +\left[\!\!
\begin{array}{c;{3pt/3pt}c}
  &  \sum\limits_{\gamma=1}^{n-1} \big(\!-\frac{\partial \psi_1^{11}}{\partial \xi_\gamma} \big)\big(\frac{\partial }{\partial x_\gamma}(\mu^{-1}g^{1\alpha}) \xi_\alpha \big)  \\
     \left[\,0\,\right]_{n\times n}& \vdots\\
    & \sum\limits_{\gamma=1}^{n-1}\big(\!- \frac{\partial \psi_1^{n-1,n-1}}{\partial \xi_\gamma} \big)\big(\frac{\partial }{\partial x_\gamma}(\mu^{-1}g^{n-1,\alpha}) \xi_\alpha\big)\\
     & \sum\limits_{\gamma=1}^{n-1}\big(\!-i\frac{\partial \psi_1^{nn}}{\partial \xi_\gamma}\big) \big( \frac{\partial }{\partial x_\gamma}(\mu^{-1}q_1^{n+1,n+1})\big) \end{array}\!\!
\right]\\ [2mm]
&& \;\;\;\,+\left[\!\!
\begin{array}{cccc;{3pt/3pt}c}
   m_1^{11} & \cdots & m_1^{1,n-1} & m_1^{1n} & m_1^{1,n+1}    \\
     \cdots & \cdots & \cdots & \cdots& \cdots \\
     m_1^{n-1,1} & \cdots & m_1^{n-1,n-1} & m_1^{n-1,n} & m_1^{n-1,n+1}  \\
     m_1^{n1} & \cdots & m_1^{n,n-1} & m_1^{nn} & m_1^{n,n+1}   \end{array}\!\!
\right].
  \end{eqnarray*}
Therefore
\begin{eqnarray}\label{2022.6.12_1} && \begin{bmatrix} \psi_0^{11}& \cdots & \psi_0^{1n}\\
\cdots & \cdots & \cdots\\
\psi_0^{n1} &\cdots & \psi_0^{nn}\!\end{bmatrix} \left[\!\!
\begin{array}{c;{3pt/3pt}c}
  &   i\mu^{-1}g^{1\alpha} \xi_\alpha   \\
     \left[\,0\,\right]_{n\times n}& \vdots\\
    & i\mu^{-1}g^{n-1,\alpha} \xi_\alpha\\
     & -\mu^{-1} q_1^{n+1,n+1} \end{array}\!\!
\right] =\left[\!\!
\begin{array}{c;{3pt/3pt}c}
  &   \phi_1^1   \\
     \left[\,0\,\right]_{n\times n}& \vdots\\
    & \phi_1^{n-1}\\
     & \phi_1^{n} \end{array}\!
\right],\end{eqnarray}
where \begin{eqnarray*} && \phi_1^{j} = \psi_1^{jj} g^{j\alpha} \frac{\partial \mu^{-1}}{\partial x_\alpha} -\sum\limits_{\gamma=1}^{n-1}\frac{\partial \psi_1^{jj}}{\partial \xi_\gamma} \frac{\partial}{\partial x_\gamma} \big(\mu^{-1} g^{j\alpha}\big) \xi_\alpha -\mu(\mu+\rho)^{-\frac{1}{2}}q_1^{j,n+1} - 2i g^{j\alpha} \xi_\alpha q_0^{n+1,n+1} \\
 && \quad\quad  \;\; - g^{j\alpha} \frac{\partial q_1^{n+1,n+1}}{\partial x_\alpha} + i \frac{\partial g^{j\alpha}}{\partial x_n} \xi_\alpha,\;\; 1\le j\le n-1;\\
 &&  \phi_1^n = \psi_1^{nn} \big( \frac{\partial \mu^{-1}}{\partial x_n} +\mu^{-1} q_0^{n+1,n+1} \big) -i \sum\limits_{\gamma=1}^{n-1}\frac{\partial \psi_1^{nn}}{\partial \xi_\gamma} \,\frac{\partial}{\partial x_\gamma} \big(\mu^{-1} q_1^{n+1,n+1}\big)
 -2\mu (\mu+\rho)^{-\frac{1}{2}}q_1^{n,n+1} \\
 && \quad \quad\;\;  + 3\sum_{l=1}^n q_0^{n+1,l} q_1^{l, n+1} -3\frac{\partial q_1^{n+1,n+1}}{\partial x_n}+6 q_1^{n+1,n+1} q_0^{n+1,n+1}
     -3i\sum\limits_{\gamma=1}^{n-1} \frac{\partial q_1^{n+1,n+1} }{\partial \xi_\gamma}\, \frac{\partial q_1^{n+1,n+1}}{\partial x_\gamma}
    \\
 && \quad \quad\;\; - \Gamma^{\beta}_{n\beta}  q_1^{n+1,n+1}+ \frac{{i}}{\sqrt{|g|}}\frac{\partial}{\partial x_\alpha}\big(\sqrt{|g|}g^{\alpha\beta}\big)\xi_\beta.\end{eqnarray*}
Here we have used the fact that $q_1^{jk}=0$ when $1\le j\ne k\le n$, and $q_1^{n+1,k}=0$ when $k\ne n+1$.
By studying above equation (\ref{2022.6.12_1}), we get that a solution of the following form:
\begin{eqnarray*}\left\{\!\!\begin{array}{ll}  \psi_0^{jk} = \frac{-i\mu \xi_k\phi_1^j}{g^{\eta\theta}\xi_\eta\xi_\theta},\,& 1\le j\le n, \; \; 1\le k\le n, \\
 \psi_0^{jn} = 0,\,& 1\le j\le n,\end{array}\right.\end{eqnarray*}
i.e.,
\begin{eqnarray}\label{2022.6.20-26}  \begin{bmatrix} \psi_0^{11}& \cdots & \psi_0^{1n}\\
\cdots & \cdots & \cdots\\
\psi_0^{n1} &\cdots &\cdots \psi_0^{nn}\!\end{bmatrix} =  \begin{bmatrix} -\frac{i\mu \xi_1 \phi_1^1}{g^{\eta\theta} \xi_\eta\xi_\theta} & \cdots & -\frac{i\mu \xi_{n-1} \phi_1^1}{g^{\eta\theta} \xi_\eta\xi_\theta} & 0\\
\cdots & \cdots & \cdots\\
-\frac{i\mu \xi_1 \phi_1^{n-1}}{g^{\eta\theta} \xi_\eta\xi_\theta} & \cdots & -\frac{i\mu \xi_{n-1} \phi_1^{n-1}}{g^{\eta\theta} \xi_\eta\xi_\theta} & 0\\
-\frac{i\mu \xi_1 \phi_1^{n}}{g^{\eta\theta} \xi_\eta\xi_\theta} & \cdots & -\frac{i\mu \xi_{n-1} \phi_1^{n}}{g^{\eta\theta} \xi_\eta\xi_\theta} & 0
\!\end{bmatrix}.\end{eqnarray}

Generally, for $m=2,3,4,\cdots$, by considering the symbol equations with homogeneous of degree $2-m$ in $\xi'$ for the $\boldsymbol{\Psi} \mathbf{L}=\mathbf{M}$, we have
\begin{eqnarray} \label{2022.7.5-2}  \;\;\quad \quad \;\;\boldsymbol{\psi}_{1-m} \left[\!\!
\begin{array}{c;{3pt/3pt}c}
  &   {i}\mu^{-1}g^{1\alpha} \xi_\alpha   \\
     \left[\,0\,\right]_{n\times n}& \vdots\\
    &{i}\mu^{-1}g^{n-1,\alpha} \xi_\alpha\\
     & -\mu^{-1} q_1^{n+1,n+1} \end{array}\!\!
\right] =- \sum_{\substack{j-|\vartheta|+k =2-m\\ (j,k,|\vartheta|)\ne (1-m,1,0)}} \frac{(-i)^{|\vartheta|}}{\vartheta!} \partial^\vartheta_{\xi'} \boldsymbol{\psi}_j \partial_{x'}^\vartheta \boldsymbol{\it l}_k +\mathbf{m}_{2-m},\end{eqnarray}
 where $\mathbf{m}_{2-m}$ is homogeneous of
 degree $2-m$ in $\xi'$ for the pseudodifferential operator $\mathbf{M}$ in (\ref{2022.7.5-1}).  Denote by $(\phi^1_{2-m}, \cdots, \phi^n_{2-m})^T$ the $(n+1)$-column of the matrix on the right-hand side of (\ref{2022.7.5-2}), we get
 \begin{eqnarray}\label{2022.6.20-6'}   \boldsymbol{\psi}_{1-m} = - \begin{bmatrix} \frac{i\mu \xi_1 \phi_{2-m}^1}{g^{\eta\theta} \xi_\eta\xi_\theta} & \cdots & \frac{i\mu \xi_{n-1} \phi_{2-m}^1}{g^{\eta\theta} \xi_\eta\xi_\theta} & 0\\
\cdots & \cdots & \cdots\\
\frac{i\mu \xi_1 \phi_{2-m}^{n-1}}{g^{\eta\theta} \xi_\eta\xi_\theta} & \cdots & \frac{i\mu \xi_{n-1} \phi_{2-m}^{n-1}}{g^{\eta\theta} \xi_\eta\xi_\theta} & 0\\
\frac{i\mu \xi_1 \phi_{2-m}^{n}}{g^{\eta\theta} \xi_\eta\xi_\theta} & \cdots & \frac{i\mu \xi_{n-1} \phi_{2-m}^{n}}{g^{\eta\theta} \xi_\eta\xi_\theta} & 0
\!\end{bmatrix}, \;\;\;\; m=1,2,3,\cdots.\end{eqnarray}
Since $\sigma_\mu (\mathbf{u},p)(-\boldsymbol{\nu})=-\sigma_\mu (\mathbf{u},p)\boldsymbol{\nu}$ on $\partial \Omega$,  the required conclusions for $\sigma_\mu (\mathbf{u},p)\boldsymbol{\nu}$ can be immediately  obtained.

\vskip 0.28 true cm

\noindent {\bf Remark 4.2.}  The above lemma implies that the pseudodifferential operator $\boldsymbol{\Lambda}$  has been explicitly obtained on $\partial \Omega$,
and the following representation holds: \begin{eqnarray*}\boldsymbol{\Lambda}  \mathbf{u_0}(x')=\int_{\mathbb{R}^{n-1}}\!\! e^{{i}\langle x',\xi'\rangle} \mathbf{\boldsymbol{\psi}}
 (x',\xi)  \hat{ \mathbf{u}}_0(\xi')\, d\xi' \;\;\;\mbox{for any}\;\, \mathbf{u}_0\in [C_0^\infty (T_{x'}(\partial \Omega))]^n \;\, \mbox{and}\;\,  x'\in \partial \Omega.\end{eqnarray*}

\vskip 0.28 true cm

\noindent {\bf Remark 4.3.}  From (\ref{2022.7.1-8}) and (\ref{2022.7.5-4}) we can see that the Dirichlet-to-Neumann map $\boldsymbol{\Lambda}$ of a Stokes flow uniquely determines the $\mu\big|_{\partial \Omega}$ and $\frac{\partial \mu}{\partial \nu}\big|_{\partial \Omega}$.

\vskip 1.48 true cm

\section{Asymptotic expansion of the heat trace of the Dirichlet-to-Neumann map}

\vskip 0.56 true cm

\noindent  {\it Proof of Theorem 1.1.} \  Since the principal symbol $2\sqrt{g^{\alpha\beta}\xi_\alpha\xi_\beta}\;\mathbf{I}_n$ of $\boldsymbol{\Lambda}$ is  positive definite, we get that $\boldsymbol{\Lambda}$ is an elliptic pseudodifferential operator of order $1$.  Let $\epsilon >0$ be given.
 It follows that the spectrum of $\boldsymbol{\Lambda}$ lies in a cone of slope $\epsilon$ about the positive real axis. Let
 $\mathcal{C}$ be a path about the cone with slope $2\epsilon$ outside some compact set. For $\tau$
on $\mathcal{C}$, the operator $(\boldsymbol{\Lambda}-\tau \mathbf{I})^{-1}$ is a uniformly bounded compact operator from
 $[L^2(\partial \Omega)]^n\to [L^2(\partial \Omega)]^n$. The integral
  \begin{eqnarray} \label{3>1} \frac{i}{2\pi} \int_{\mathcal{C}} e^{-t\tau } (\boldsymbol{\Lambda} - \tau \mathbf{I})^{-1} d\tau\end{eqnarray}
converges absolutely for $t>0$ and defines the operator semigroup  $e^{-t\boldsymbol{\Lambda}}$ (see \cite{Gr} and  \cite{Pa}).
We construct a matrix-valued pseudodifferential operator to approximate the resolvant $(\boldsymbol{\Lambda} -\tau \mathbf{I})^{-1}$ as follows: let  $\mathbf{\mathcal{S}}(\tau)$ be a matrix-valued pseudodifferential operator of order $-1$ with parameter $\tau$ for which
\begin{eqnarray*} (\boldsymbol{\Lambda}-\tau \mathbf{I}) \mathbf{\mathcal{S}}(\tau) =\mathbf{\mathcal{S}}(\tau)(\boldsymbol{\Lambda}-\tau \mathbf{I}) =\mathbf{I}.\end{eqnarray*}
(Actually, we require that $(\boldsymbol{\Lambda}-\tau \mathbf{I}) \mathbf{\mathcal{S}}(\tau)-\mathbf{I}$ and $\mathbf{\mathcal{S}}(\tau)(\boldsymbol{\Lambda}-\tau \mathbf{I}) -\mathbf{I}$ are both pseudodifferential operator of order $-\infty$.)  Let the  pseudodifferential operator $\mathbf{\mathcal{S}}(\tau)$  has the symbol expansion
\begin{eqnarray*} \iota(\mathbf{\mathcal{S}}(\tau))\sim \sum_{l=0}^\infty \boldsymbol{\varpi}_{-1-l} (x', \xi', \tau).\end{eqnarray*}
 Then
\begin{eqnarray*} \iota((\boldsymbol{\Lambda}-\tau \mathbf{I}) \mathbf{\mathcal{S}}(\tau))=\sum_{\vartheta}\frac{(-i)^{|\vartheta|}}{\vartheta!} \big(\partial_{\xi'}^\vartheta \iota(\boldsymbol{\Lambda}-\tau \mathbf{I})\big)(\partial_{x'}^\vartheta \iota(\mathbf{\mathcal{S}}(\tau))).\end{eqnarray*}
Since  $\iota(\boldsymbol{\Lambda})\sim \sum_{m=0}^\infty \boldsymbol{\psi}_{1-m}$, we can decompose the above sum into order of
homogeneity  \begin{eqnarray*}\iota((\boldsymbol{\Lambda}-\tau \mathbf{I}) \mathbf{\mathcal{S}(\tau)})\sim  \sum\limits_{\substack{-l=k-|\vartheta| -1-m\\ 0\le m\le l, \, k\le 1,\,\,-1-l \le -1}}
 \frac{(-i)^{|\vartheta|}}{\vartheta!} \, \big(\partial^\vartheta_{\xi'} \tilde{\boldsymbol{\psi}}_k\big)(\partial_{x'}^\vartheta \boldsymbol{\varpi}_{-1-m}), \end{eqnarray*}
  where \begin{eqnarray*} \tilde{\boldsymbol{\psi}}_k(x', \xi',\tau)=\left\{ \begin{array}{ll} \boldsymbol{\psi}_1 -\tau \mathbf{I}_n\,\;&\mbox{when}\;\, k=1,\\
 \boldsymbol{\psi}_k \;\; &\mbox{when}\;\; k=0, -1, -2, \cdots.\end{array}\right.\end{eqnarray*}
These equations define the $\boldsymbol{\varpi}_{-1-l}$ inductively, as follows:
 \begin{eqnarray*}
 & \boldsymbol{\varpi}_{-1} (x', \xi', \tau) = (\boldsymbol{\psi}_1 -\tau \mathbf{I}_n)^{-1},
\\ &  \boldsymbol{\varpi}_{-1-l} (x', \xi', \tau) =- (\boldsymbol{\psi}_1- \tau \mathbf{I}_n)^{-1} \left(\sum\limits_{\substack{-l=k-|\vartheta| -1-m\\0\le m<l, \, k\le 1,\,\,l\ge 1}}
 \frac{(-i)^{|\vartheta|}}{\vartheta!} \, (\partial^\vartheta_{\xi'} \tilde{\boldsymbol{\psi}}_k)(\partial_{x'}^\vartheta \boldsymbol{\varpi}_{-1-m})\right), \;\, l\ge 1.\end{eqnarray*}
 For the sake of convenience, we also write out the expression for $\boldsymbol{\varpi}_{-2}$:
 \begin{eqnarray} \label{2022.6.27-13} &&\boldsymbol{\varpi}_{-2}(x', \xi', \tau) = - ({\boldsymbol{\psi}}_1 -\tau \mathbf{I}_n)^{-1} \boldsymbol{\psi}_0 \boldsymbol{\varpi}_{-1}-i (\boldsymbol{\psi}_1-\tau \mathbf{I}_n)^{-1} \sum_{|\vartheta|=1} \big(\partial^\vartheta_{\xi'} \boldsymbol{\psi}_1 \big) \big(\partial_{x'}^\vartheta \boldsymbol{\varpi}_{-1}\big)\\
 &&\qquad \, = - ({\boldsymbol{\psi}}_1 -\tau \mathbf{I}_n)^{-1} \boldsymbol{\psi}_0 ({\boldsymbol{\psi}}_1 -\tau \mathbf{I}_n)^{-1}-i (\boldsymbol{\psi}_1-\tau \mathbf{I}_n)^{-1} \sum_{|\vartheta|=1} (\partial^\vartheta_{\xi'} \boldsymbol{\psi}_1 )\big( \partial_{x'}^\vartheta (({\boldsymbol{\psi}}_1 -\tau \mathbf{I}_n)^{-1})\big).\nonumber\end{eqnarray}
 From (\ref{2022.6.10-5}) we see \begin{align*} & \boldsymbol{\psi}_1- \tau \mathbf{I}_n= \Big(2\sqrt{g^{\alpha \beta} \xi_\alpha \xi_\beta}-\tau\Big) \mathbf{I}_n,\end{align*}
 so that
 \begin{eqnarray} \label{2022.6.20-3}  (\boldsymbol{\psi}_1- \tau \mathbf{I}_n)^{-1}=  \frac{1}{2\sqrt{g^{\alpha \beta} \xi_\alpha \xi_\beta} -\tau} \,\mathbf{I}_n. \end{eqnarray}
It follows that
\begin{eqnarray} \label{2022.6.20-9} \mbox{Tr} \big(
 \boldsymbol{\varpi}_{-1} (x', \xi', \tau)\big)=
\mbox{Tr} \big((\boldsymbol{\psi}_1- \tau \mathbf{I}_n)^{-1}\big)=\frac{n}{
2\sqrt{g^{\alpha \beta} \xi_\alpha \xi_\beta} -\tau}.\end{eqnarray}

Next, we will calculate $\mbox{Tr} (\boldsymbol{\varpi}_{-2}(x',\xi', \tau))$.
From  (\ref{2022.7.5-4}) and (\ref{2022.6.20-3}), we have
\begin{eqnarray*}&& \!\!\!\!\!\!\!\! \!\! (\boldsymbol{\psi}_1-\tau \mathbf{I})^{-1} \boldsymbol{\psi}_0 (\boldsymbol{\psi}_1-\tau \mathbf{I})^{-1} =  \frac{\mathbf{I}_n}{2\sqrt{g^{\alpha \beta} \xi_\alpha \xi_\beta} -\tau}
  \begin{bmatrix} \frac{i\mu \xi_1 \phi_1^1}{g^{\eta\theta} \xi_\eta\xi_\theta} & \cdots & \frac{i\mu \xi_{n-1} \phi_1^1}{g^{\eta\theta} \xi_\eta\xi_\theta} & 0\\
\cdots & \cdots & \cdots\\
\frac{i\mu \xi_1 \phi_1^{n-1}}{g^{\eta\theta} \xi_\eta\xi_\theta} & \cdots & \frac{i\mu \xi_{n-1} \phi_1^{n-1}}{g^{\eta\theta} \xi_\eta\xi_\theta} & 0\\
\frac{i\mu \xi_1 \phi_1^{n}}{g^{\eta\theta} \xi_\eta\xi_\theta} & \cdots & \frac{i\mu \xi_{n-1} \phi_1^{n}}{g^{\eta\theta} \xi_\eta\xi_\theta} & 0
\!\end{bmatrix}
   \frac{\mathbf{I}_n}{2\sqrt{g^{\alpha \beta} \xi_\alpha \xi_\beta} -\tau} \\ [1mm]
  &&\qquad \qquad \qquad \qquad\;\; \qquad = \begin{bmatrix} \frac{i\mu \xi_1 \phi_1^{1}}{(2\sqrt{g^{\alpha\beta} \xi_\alpha\xi_\beta} -\tau)^2 (g^{\eta\theta} \xi_\eta\xi_\theta)} &\cdots&
  \frac{i\mu \xi_{n-1} \phi_1^{1}}{(2\sqrt{g^{\alpha\beta} \xi_\alpha\xi_\beta} -\tau)^2 (g^{\eta\theta} \xi_\eta\xi_\theta)} &0 \\
  \vdots &  \vdots & \vdots & \vdots\\
   \frac{i\mu \xi_{1} \phi_1^{n-1}}{(2\sqrt{g^{\alpha\beta} \xi_\alpha\xi_\beta} -\tau)^2 (g^{\eta\theta} \xi_\eta\xi_\theta)} &\cdots&
  \frac{i\mu \xi_{n-1} \phi_1^{n-1}}{(2\sqrt{g^{\alpha\beta} \xi_\alpha\xi_\beta} -\tau)^2 (g^{\eta\theta} \xi_\eta\xi_\theta)}&0\\
  0 & \cdots & 0 &0
  \end{bmatrix},
  \end{eqnarray*}
  so that
\begin{eqnarray} \label{2022.6.20-8} \mbox{Tr} ( (\boldsymbol{\psi}_1-\tau \mathbf{I})^{-1} \psi_0 (\boldsymbol{\psi}_1-\tau \mathbf{I})^{-1}\big)= \frac{i\mu ( \xi_1\phi_1^{1} +\cdots +\xi_{n-1} \phi_1^{n-1})}{
\big(2\sqrt{g^{\alpha \beta}\xi_\alpha\xi_\beta}-\tau \big)^2( g^{\eta\theta}\xi_\eta\xi_\theta)}.\end{eqnarray}
 By a direct calculation, we get
\begin{eqnarray} \label{2022.6.27-11}\end{eqnarray}
\begin{eqnarray*}&& \frac{i\mu(\xi_1\phi_1^1 +\cdots +\xi_{n-1}\phi_1^{n-1} )}{g^{\eta\theta}\xi_\eta\xi_\theta} =
 \frac{1}{g^{\eta\theta} \xi_\eta \xi_\theta} \left[ 2i\mu \sqrt{g^{\sigma\delta}\xi_\sigma\xi_\delta} \; g^{\beta \alpha}\xi_\beta \frac{\partial \mu^{-1}}{\partial x_\alpha}\right.\\
  && \left.\; \;\;-2i \mu \sum_\gamma \frac{\partial \sqrt{g^{\sigma\delta} \xi_\sigma\xi_\delta}}{\partial \xi_\gamma} \, \frac{\partial(\mu^{-1} g^{\beta\alpha} )}{\partial x_\gamma}\xi_\alpha \xi_\beta -i \mu^2 (\mu+\rho)^{-\frac{1}{2}} \big(q_1^{\gamma, n+1} \xi_\gamma \big)+ 2\mu g^{\beta\alpha} \xi_\alpha \xi_\beta q_0^{n+1,n+1}\right.\\
  &&  \left. \;\; \; -i \mu g^{\alpha\beta} \xi_\beta \frac{\partial q_1^{n+1,n+1}}{\partial \xi_\alpha} - \mu \frac{\partial g^{\beta\alpha}}{\partial x_n} \xi_\beta \xi_\alpha\right] \\
&& =2\mu q_0^{n+1,n+1} - \frac{\mu \frac{\partial g^{\alpha\beta}}{\partial x_n} \xi_\alpha\xi_\beta}{g^{\eta\theta}\xi_\eta \xi_\theta}
 -i\mu \frac{g^{\alpha\beta}\xi_\beta}{g^{\eta\theta}\xi_\eta\xi_\theta} \frac{\partial q_1^{n+1,n+1}}{\partial x_\alpha} - \frac{2i }{g^{\eta\theta}\xi_\eta\xi_\theta} \sum_\gamma \frac{\partial \sqrt{g^{\sigma\delta}\xi_\sigma\xi_\delta}}{\partial \xi_\gamma}\, \frac{\partial g^{\alpha\beta}}{\partial x_\gamma}\xi_\alpha\xi_\beta\\
 &&\quad \;- i\mu^2 (\mu+\rho)^{-\frac{1}{2}} \frac{q_1^{\gamma,n+1}\xi_\gamma}{g^{\eta\theta}\xi_\eta\xi_\theta}.  \end{eqnarray*}

Since the $(n+1)$-th row of $\mathbf{A}_1$ vanishes, we get that the $(n+1)$-th row of $\mathbf{A}_1\mathbf{E}_1$ and $\mathbf{A}_1\mathbf{E}_1\mathbf{A}_2$ also vanish. Note that
the $(n+1)$-th row of $\mathbf{E}_1$ is \begin{align*} & \Big(\! \!-i\mu (\mu+\rho)^{-\frac{1}{2}}(\xi_1,\,\cdots,\, \xi_{n-1}), \; \mu (\mu+\rho)^{-\frac{1}{2}} \sqrt{g^{\eta\theta} \xi_\eta\xi_\theta}, \; i\sum_{\gamma} \frac{\partial q_1^{n+1,n+1}}{\partial \xi_\gamma} \frac{\partial q_1^{n+1,n+1}}{\partial x_\gamma} + \frac{\partial \sqrt{g^{\eta\theta}\xi_\eta\xi_\theta}}{\partial x_n}\\
&+\Gamma_{n\beta}^\beta \sqrt{g^{\eta\theta} \xi_\eta\xi_\theta} + \frac{\Gamma_{\gamma \sigma}^n g^{\alpha\gamma}g^{\beta \sigma} \xi_\alpha \xi_\beta}{\sqrt{g^{\eta\theta}\xi_\eta\xi_\theta}} - i\big(g^{\alpha\beta} \Gamma_{\alpha\gamma}^\gamma +\frac{\partial g^{\alpha\beta}}{\partial x_\alpha}\big) \xi_\beta\; \Big)\end{align*}
and  \begin{eqnarray*}
  \mathbf{A}_2\!\! &=\!\!&\begin{small}\begin{bmatrix}
\begin{BMAT}(@, 15pt, 15pt){c.c}{c.c}
   \!\left[0\right]_{n\times n}   &   \left[2{i}\,  \Gamma_{\beta n}^j g^{\beta \alpha}\xi_\alpha\!+\!
   \frac{1}{\sqrt{ g^{\eta\theta}\xi_\eta\xi_\theta }}\Gamma_{\gamma \sigma}^j g^{\alpha \gamma} g^{\beta \sigma} \xi_\alpha\xi_\beta \right]_{n\times 1} \\
\left[0 \right]_{1\times n}   &  0  \end{BMAT}
\end{bmatrix}.\end{small}\end{eqnarray*}
Therefore, the $(n+1)$-th row of $\mathbf{E}_1\mathbf{A}_2$ is
$$\Big( \underset{n}{\underbrace{
0,\cdots,0}}, \mu (\mu+\rho)^{-\frac{1}{2}} \Gamma^{n}_{\gamma \sigma}g^{\alpha\gamma} g^{\beta\sigma}\xi_\alpha\xi_\beta +2 \mu (\mu+\rho)^{-\frac{1}{2}} \xi_\gamma \Gamma^\gamma_{\beta n} g^{\beta \alpha}\xi_\alpha - \frac{i \mu (\mu+\rho)^{-\frac{1}{2}}\xi_\delta}{\sqrt{g^{\eta\theta}\xi_\eta\xi_\theta}} \Gamma^\delta_{\gamma \sigma} g^{\alpha\gamma}g^{\beta\sigma}\xi_\alpha\xi_\beta\Big).$$
From (\ref{2022.6.23-5}) we see that the $(n+1,n+1)$ entry of $\mathbf{q}_0$ is
\begin{eqnarray*} &&\!\!\!\!\! \!\!\!\!\!q_0^{n+1,n+1}\!=\! \frac{1}{2\sqrt{g^{\eta\theta}\xi_\eta\xi_\theta}} \bigg( i\sum_\gamma\frac{\partial q_1^{n+1,n+1}}{\partial \xi_\gamma}\frac{\partial q_1^{n+1,n+1}}{\partial x_\gamma}\!
+\! \frac{\partial \sqrt{g^{\sigma\delta}\xi_\sigma\xi_\delta}}{\partial x_n} +\Gamma_{n\beta}^\beta \sqrt{g^{\sigma\delta}\xi_\sigma\xi_\delta} +
\frac{\Gamma_{\gamma\sigma}^n g^{\alpha \gamma}g^{\beta \sigma} \xi_\alpha \xi_\beta} {\sqrt{g^{\delta\varepsilon}\xi_\delta\xi_\varepsilon}} \\
&& \qquad \quad -i\big( g^{\alpha\beta} \Gamma_{\alpha\gamma} ^\gamma +\frac{\partial g^{\alpha\beta}}{\partial x_\alpha}\big) \xi_\beta  \bigg) \!-\! \frac{(\mu+\rho)^{\frac{1}{2}}}{4\mu g^{\eta\theta} \xi_\eta\xi_\theta} \bigg( \frac{\mu}{(\mu+\rho)^{\frac{1}{2}}} \Gamma^n_{\gamma \sigma} g^{\alpha\gamma} g^{\beta \sigma} \xi_\alpha\xi_\beta \\
 && \qquad  \quad +2 \mu(\mu+\rho)^{-\frac{1}{2}} \Gamma^\gamma_{\beta n} g^{\beta \alpha} \xi_\gamma \xi_\alpha - \frac{i\mu (\mu+\rho)^{-\frac{1}{2}} }{\sqrt{g^{\eta\theta}\xi_\eta\xi_\theta}} \Gamma^{\delta}_{\gamma\sigma} g^{\alpha\gamma} g^{\beta\sigma} \xi_\delta \xi_\alpha \xi_\beta\bigg).\end{eqnarray*}
 Using the equality $\frac{\partial g_{\gamma \sigma}}{\partial x_n} g^{\sigma \beta}=- g_{\gamma \sigma} \frac{\partial g^{\sigma\beta} }{\partial x_n}$, we get
 \begin{align*} && q_0^{n+1,n+1}  = \frac{ i}{2\sqrt{g^{\eta\theta}\xi_\eta \xi_\theta} }\sum_\gamma \frac{\partial q_1^{n+1,n+1}}{\partial \xi_\gamma}\frac{\partial q_1^{n+1,n+1}}{\partial x_\gamma} -\frac{1}{4} g_{\beta \gamma}\frac{\partial g^{\beta \gamma}}{\partial x_n} +\frac {5}{8 g^{\eta\theta}\xi_\eta\xi_\theta}  \, \frac{\partial g^{\alpha\beta}}{\partial x_n} \xi_\alpha \xi_\beta \\
&& \qquad \qquad  -
 \frac{1}{2\sqrt{g^{\eta\theta}\xi_\eta \xi_\theta}} \big( g^{\alpha\beta} \Gamma_{\alpha\gamma}^\gamma +\frac{\partial g^{\alpha\beta}}{\partial x_\alpha} \big) \xi_\beta +\frac{i}{4(g^{\eta\theta} \xi_\eta\xi_\theta)^{3/2}} \xi_\delta \Gamma_{\gamma \sigma}^\delta g^{\alpha \gamma} g^{\beta \sigma} \xi_\alpha\xi_\beta.\end{align*}
  At each point $x_0=(x'_0,0)\in \partial \Omega$, in the boundary normal coordinates  we  have $\frac{\partial g^{\alpha\beta}}{\partial x_\gamma}(x_0)=0$, $g^{n\alpha}=0$, so that
 \begin{eqnarray*}  q_0^{n+1, n+1}(x_0)\!\!\!\!\!\!\!\!&&=
  -\frac{1}{4} \sum_\alpha \frac{\partial g^{\alpha \alpha}}{\partial x_n} (x_0) + \frac{5}{8|\xi'|^2} \frac{\partial g^{\alpha \beta}}{\partial x_n} (x_0) \xi_\alpha\xi_\beta .\end{eqnarray*}
Hence
  \begin{eqnarray*} \frac{i\mu(\xi_1\phi_1^1 +\cdots +\xi_{n-1}\phi_1^{n-1} )}{g^{\eta\theta}\xi_\eta\xi_\theta}\Big|_{x=x_0} =
 2 \mu q_0^{n+1,n+1} (x_0)
 -\frac{\mu \frac{\partial g^{\alpha\beta}}{\partial x_n}(x_0)\, \xi_\alpha\xi_\beta}{g^{\eta \theta}(x_0) \xi_\eta\xi_\theta }-i \mu^2 (\mu+\rho)^{-\frac{1}{2}} \frac{\big(q_1^{\gamma, n+1}(x_0)\, \xi_\gamma \big)}{g^{\eta \theta}(x_0) \xi_\eta\xi_\theta} .\end{eqnarray*}
 Similarly, we have  $\partial_{x'}^\alpha \boldsymbol{\varpi}_{-1}(x_0) =\partial_{x'}^\alpha (\boldsymbol{\psi}_1-\tau \mathbf{I})^{-1} \big|_{x=x_0} =0$, so that
 $\mbox{Tr}\, \big((\boldsymbol{\psi}_1-\tau I_n)^{-1} \sum_{|\alpha|=1} \partial^\alpha_{\xi'} \boldsymbol{\psi}_1 \cdot D_{x'}^\alpha \boldsymbol{\varpi}_{-1}\big)=0$ at  $x_0\in \partial \Omega$.
  That is, the trace of the second term in (\ref{2022.6.27-13}) also vanishes.
  Note that \begin{eqnarray*} -i\mu^2 (\mu+\rho)^{-\frac{1}{2}} \frac{q_1^{\gamma,n+1}(x_0)\xi_\gamma}{|\xi'|^2} = \frac{1}{|\xi'|^2} \Big( 2\mu \Gamma_{\beta n}^\delta g^{\beta \alpha} \xi_\alpha \xi_\delta -\frac{i \mu \Gamma_{\gamma \sigma}^\delta g^{\alpha \gamma}g^{\beta \sigma} \xi_\alpha \xi_\beta\xi_\delta}{\sqrt{g^{\eta\theta}\xi_\eta \xi_\theta}}\Big) \Big|_{x=x_0} = - \frac{\mu\, \frac{\partial g^{\alpha\beta}}{\partial x_n} (x_0) \xi_\alpha \xi_\beta}{|\xi'|^2}.\end{eqnarray*}
  Combining these, (\ref{2022.6.27-13}), (\ref{2022.6.20-8}), (\ref{2022.6.27-11}) and (\ref{18/7/14/1}) we get that at each $x_0\in \partial \Omega$,
 \begin{eqnarray} \label{2022.6.27-12}&& \;\;\;\;  \boldsymbol{\varpi}_{-2}(x', \xi',\tau)=\mbox{Tr} ( (\boldsymbol{\psi}_1-\tau \mathbf{I})^{-1} \psi_0 (\boldsymbol{\psi}_1-\tau \mathbf{I})^{-1}\big)\\
&&\; \;= \frac{1}{(2|\xi'|^2-\tau)^2}\bigg[
-\frac{\mu}{2} \sum_{\alpha=1}^{n-1} \frac{\partial g^{\alpha\alpha}}{\partial x_n} (x_0) -  \frac{3\mu \frac{\partial g^{\alpha\beta}}{\partial x_n}(x_0) \xi_\alpha\xi_\beta}{4|\xi'|^2} \bigg],\nonumber \end{eqnarray}
where $|\xi'|= \sqrt{\sum_{\alpha=1}^{n-1} \xi_\alpha^2}$.

Furthermore,
   there exists a matrix-valued function $\mathbf{K} (t, x', y')$, which is called the parabolic (or semigroup) kernel, such that (see \cite{Gr}, or p.$\,$4 of \cite{Frie})
                \begin{eqnarray*}  e^{-t\boldsymbol{\Lambda}}\,\mathbf{u}_0(x')=\int_{\partial \Omega} \mathbf{K}(t, x',y') \mathbf{u}_0(y')\,ds(y'), \quad \,
        {\mathbf{w}}_0\in [H^{\frac{1}{2}}(\partial \Omega)]^n, \end{eqnarray*}
 Recall that $\{{\boldsymbol{\phi}}_k\}_{k=1}^\infty$ be orthnormal eigenvectors of the Dirichlet-to-Neumann map $\Xi_g$ corresponding to eigenvalues $\{\lambda_k\}_{k=1}^\infty$, then the parabolic kernel  ${\mathbf{K}}(t, x', y')=e^{-t \boldsymbol{\Lambda}} \delta(x'-y')$ is given by \begin{eqnarray} \label{18/12/18} {\mathbf{K}}(t,x',y') =\sum_{k=1}^\infty e^{-t \lambda_k} {\boldsymbol{\phi}}_k(x')\otimes {\boldsymbol{\phi}}_k(y').\end{eqnarray}
This implies that the integral of the trace of ${\mathbf{K}}(t,x',y')$ is actually a spectral invariants:
\begin{eqnarray} \label{1-0a-2}\int_{\partial\Omega} \mbox{Tr}({\mathbf{K}}(t,x',x'))\, ds(x')=\sum_{k=1}^\infty e^{-t \lambda_k}.\end{eqnarray}
 It follows that \begin{eqnarray}\label{18/12/2} {\mathbf{K}}  (t,x',y') \!\!&\!=\!&\! e^{-t\boldsymbol{\Lambda}}\big(\delta(x'-y')\mathbf{I}_{n-1}\big)  \\
&\!=\!&\!\!\! \frac{1}{(2\pi)^{n-1}} \!\int_{T^*(\partial \Omega)} \!\!e^{i\langle x'-y', \xi'\rangle} \bigg\{\!\frac{i}{2\pi}\! \int_{\mathcal{C}} \!e^{-t\tau}\! \big(\iota((\boldsymbol{\Lambda}\!-\!\tau \mathbf{I})^{-1})\big) d\tau\!\bigg\} d\xi'\nonumber\\
 &\!=\!&\!\! \! \frac{1}{(2\pi)^{n-1}} \int_{\mathbb{R}^{n-1}} e^{i\langle x'-y', \xi'\rangle} \bigg\{\frac{i}{2\pi} \int_{\mathcal{C}} e^{-t\tau} \begin{small} \bigg(\end{small}
 \boldsymbol{\varpi}_{-1} (x',\xi',\tau) \nonumber\\
 \!\!&&\! \!+\boldsymbol{\varpi}_{-2} (x',\xi',\tau) +\boldsymbol{\varpi}_{-3} (x',\xi',\tau) +\cdots \begin{small}\bigg) \end{small} d\tau\bigg\} d\xi',\nonumber\end{eqnarray}
 where $T^*(\partial \Omega)$ is the cotangent space at $x'$, so that \begin{eqnarray} \label{18/12/2-2}  && \mbox{Tr}\,(\mathbf{{K}}  (t,x',x')) \!\!=\!\! \frac{1}{(2\pi)^{n-1}}\!\int_{\mathbb{R}^{n-1}} \!\bigg\{\!\frac{i}{2\pi}\! \int_{\mathcal{C}}\! e^{-t\tau}\!  \bigg(\!
 \mbox{Tr}\,\big(\boldsymbol{\varpi}_{-1} (x',\xi',\tau)\big)\\
 && \quad \;\;+ \mbox{Tr}\,\big(\boldsymbol{\varpi}_{-2} (x',\xi',\tau)\big)\! +\! \mbox{Tr}\,\big(\boldsymbol{\varpi}_{-3} (x',\xi',\tau)\big) \!+\!\cdots \bigg) d\tau\!\bigg\} d\xi'.\nonumber\end{eqnarray}
 Combining (\ref{1-0a-2}) and (\ref{18/12/2-2}), we get that \begin{eqnarray} && \sum_{k=1}^\infty e^{-t \lambda_k}=\int_{\partial \Omega} \bigg\{ \frac{1}{(2\pi)^{n\!-\!1}}\!  \int_{\mathbb{R}^{n\!-\!1}}\! \bigg[\frac{i}{2\pi} \!\int_{\mathcal{C}}\! e^{-t\tau} \! \bigg(\!
 \mbox{Tr}\,\big(\boldsymbol{\varpi}_{-1} (x',\xi',\tau)\big) \\
 &&\quad \quad \;\quad \quad \quad+ \,\mbox{Tr}\,\big(\boldsymbol{\varpi}_{-2} (x',\xi',\tau)\big)\!+\! \mbox{Tr}\,\big(\boldsymbol{\varpi}_{-3} (x',\xi',\tau)\big) \!+\!\cdots \!\bigg) d\tau\!\bigg] d\xi'\bigg\} dS(x'),\nonumber\end{eqnarray}
 It follows that for $ 0\le l< n$,
 \begin{eqnarray} \label{3./1} && a_{l}(x) =\frac{i}{(2\pi)^{n}}  \int_{T^*_{x} (\partial \Omega)} \int_{\mathcal{C} }\mbox{Tr}\,\big(e^{-\tau } \boldsymbol{\varpi}_{-1-l} (x', \xi',\tau)\big) d\tau \, d\xi'\\
&& \quad \quad \;\; =\frac{1}{(2\pi)^{n-1}}  \int_{\mathbb{R}^{n-1}} \bigg(\frac{i}{2\pi} \int_\Gamma e^{-\tau }\, \mbox{Tr}\,\big(\boldsymbol{\varpi}_{-1-l} (x', \xi',\tau)\big) d\tau \bigg) d\xi'.\nonumber\end{eqnarray}
Hence
\begin{align*} \sum_{k=1}^\infty e^{-t\lambda_k} &\!=  \int_{\partial \Omega} \left\{ \frac{1}{(2\pi)^{n-1}} \!\int_{\mathbb{R}^{n-1}} \!\left[\frac{i}{2\pi} \int_{\mathcal{C}} e^{-t\tau}  \Big(  \mbox{Tr}\big((\boldsymbol{\psi}_1-\tau \mathbf{I})^{-1} \big)\right.\right.\\
 &\left.\left.\;\,\;\; + \mbox{Tr} \big( (\boldsymbol{\psi}_1-\tau \mathbf{I})^{-1} \boldsymbol{\psi}_0 (\boldsymbol{\psi}_1\!-\!\tau \mathbf{I})^{-1}\big)+\cdots \Big) d\tau\right]d\xi' \right\} ds(x')\\
&\!=   \int_{\partial \Omega}\! \left\{ \!\frac{1}{(2\pi)^{n\!-1}} \!\int_{\mathbb{R}^{n-1}}\!\left[ \frac{i}{2\pi} \int_{\mathcal{C}} e^{-t\tau} \Big( \frac{n}{2|\xi'|-\tau} \Big)\right] \!d\xi'\! \right\} ds(x') \\
&\;\; \;\; + \!\int_{\partial \Omega}\! \bigg\{\! \frac{1}{(2\pi)^{n-1}} \!\int_{\mathbb{R}^{n\!-1}}\!\bigg[\! \frac{i}{2\pi} \!\int_{\mathcal{C}} e^{-t\tau}\frac{1}{(2|\xi'|\!-\!\tau)^2}\,\Big(\!\!-\!\mu\sum_{\alpha=1}^{n-1} \kappa_\alpha(x') \\
&\;\;\;\;  -\frac{3\mu \sum_{\alpha=1}^{n-1}\kappa_\alpha (x') \xi_\alpha^2} {2|\xi'|^2}\Big) d\tau \bigg] d\xi' \bigg\} ds(x')+\cdots
 \end{align*}
Applying residue theorem (see, for example, p.$\,$147-151 of \cite{Ah}) and  the following formulas (see \cite{Liu4})
\begin{align*}& \;\int_{\mathbb{R}^{n-1}}e^{-c|\xi|}\,d\xi = \frac{\Gamma(n-1)\operatorname{vol}(\mathbb{S}^{n-2})}{c^{n-1}},\quad c>0,\ n\geqslant 2,\\
    &\int_{\mathbb{R}^{n-1}}e^{-c|\xi|} \frac{\xi_{\alpha}^2}{|\xi|^2} \,d\xi = \frac{\Gamma(n-1)\operatorname{vol}(\mathbb{S}^{n-2})}{(n-1)c^{n-1}},\quad c>0,\ n\geqslant 2,
\end{align*} we find  that, as $t\to 0^+$,
 \begin{align*}
\sum_{k=1}^\infty e^{-t\lambda_k} & = \int_{\partial \Omega}\! \bigg\{ \!\frac{1}{(2\pi)^{n\!-1}} \!\int_{\mathbb{R}^{n-1}} \Big(n e^{-2t|\xi|} \Big) d\xi' \bigg\}
  ds(x') \\
  & \quad  + \! \int_{\partial \Omega}\! \bigg\{ \!\frac{1}{(2\pi)^{n\!-1}} \!\int_{\mathbb{R}^{n-1}}\! \bigg[\!-t e^{-2t|\xi|} \Big(\!-\mu\sum_{\alpha=1}^{n-1} \kappa_\alpha(x')-\frac{3\mu \sum_{\alpha=1}^{n-1} \kappa_\alpha(x')  \xi_\alpha^2} {2|\xi'|^2} \Big) \bigg] d\xi' \bigg\}
  ds(x') +\cdots\\
       & =  \frac{n}{(2\pi)^{n-1}}\,\frac{\Gamma(n-1) \mbox{vol}(\mathbb{S}^{n-2})}{(2 t)^{n-1}}\! \int_{\partial \Omega}\! ds(x')
       +\frac{t\,\Gamma(n\!-\!1) \, \mbox{vol}(\mathbb{S}^{n-2})}{(2\pi)^{n-1} (2t)^{n-1} }\!\int_{\partial \Omega}\! \mu(x') \big(\sum_{\alpha=1}^{n-1} \kappa_\alpha(x')\big) ds(x')
       \\
    & \quad \, + \! \int_{\partial \Omega}\!\frac{3}{2(2\pi)^{n-1}} \cdot   \frac{t\,\Gamma(n-1) \mbox{vol}(\mathbb{S}^{n-2}) }{ (n-1)(2 t)^{n-1}}\, \mu \big(\sum_{\alpha=1}^{n-1} \kappa_\alpha(x') \big) \,ds(x')  + \left\{\!\begin{array}{ll} O(t^{3-n}) \quad \;\, \mbox{when}\;\; n>2,\\ O(t\log t) \quad \, \mbox{when} \;\; n=2\end{array}\right. \\
   & = \frac{n\, \Gamma(n\!-\!1) \mbox{vol}(\mathbb{S}^{n-2})}{(4\pi t)^{n-1}} \,\mbox{vol}(\partial \Omega)\! + \!  \frac{(2n+1)\,t\,\Gamma(n\!-\!1) \mbox{vol}(\mathbb{S}^{n-2}) }{2(n-1)(4\pi t)^{n-1}} \! \! \int_{\partial \Omega}\!\mu(x') \big(\sum_{\alpha=1}^{n-1} \kappa_\alpha (x')\big)\, ds(x')  \\
    & \quad  + \left\{\!\begin{array}{ll} O(t^{3-n}) \quad \;\, \mbox{when}\;\; n>2,\\ O(t\log t) \quad \, \mbox{when} \;\; n=2.\end{array}\right.
  \end{align*}    \qed

\noindent {\bf Remark 5.1.} Clearly, Theorem 1.1 holds when $\Omega\subset \mathbb{R}^n$ is a bounded domain with smooth boundary.

\vskip 1.18  true cm

\centerline {\bf  Acknowledgments}

\vskip 0.38 true cm

\vskip 0.38 true cm
The author would like to thank Professors Fang-Hua Lin, Jalal Shatah, Gunther Uhlmann, Ari Laptev,  Mark S. Ashbaugh, Robert Kohn and Jenn-Nan Wang for their great support and many useful comments and discussions. This research was supported by NNSF of China (11671033/A010802).

  \vskip 1.58 true cm

\end{document}